\documentclass[11pt,twoside]{atmp}

\newtheorem{theorem}{Theorem}[section]
\newtheorem{lemma}[theorem]{Lemma}

\theoremstyle{definition}

\theoremstyle{remark}

\numberwithin{equation}{section}

\input{tcilatex}
\begin{document}

\title[Fokker-Planck dynamics for the Ricci flow]
{Fokker-Planck dynamics and entropies for
 the normalized Ricci flow}

\author{MAURO CARFORA}
\address{Dipartimento di Fisica Nucleare e Teorica, Universita` degli Studi di Pavia\\
and\\ 
Istituto Nazionale di Fisica Nucleare, Sezione di Pavia\\
via A. Bassi 6, I-27100 Pavia, Italy}
\email{mauro.carfora@pv.infn.it}
\thanks{Research supported in part by PRIN Grant $\# 2004012375-002$.}



\date{}



\begin{abstract}
We consider some elementary aspects of the geometry of the space of probability measures endowed with Wasserstein distance. In such a setting,  we discuss the various terms entering Perelman's shrinker entropy, and characterize two new monotonic  functionals for the volume-normalized Ricci flow. One is obtained by a rescaling of the curvature term in the shrinker entropy. The second is associated with a gradient flow obtained by adding a curvature-drift  to Perelman's backward heat equation. We show that the resulting Fokker-Planck PDE  is the natural diffusion flow for probability measures absolutely continuous with respect to the Ricci-evolved Riemannian measure, we discuss its exponential trend to equilibrium and its relation with the viscous Hamilton--Jacobi equation.
\end{abstract}

\maketitle

\section*{INTRODUCTION}
The Ricci flow introduced by R. Hamilton \cite{11}, (see \cite{1,6,12} for reviews), is a geometric evolution equation which deforms the metric $g$ of a Riemannian manifold $(\Sigma ,g)$ in the direction of its Ricci curvature $Ric(g)$. Under suitable conditions it provides the natural technique for smoothing and uniformizing $(\Sigma,g)$  to specific model geometries. From the perspective of theoretical physics, the Ricci flow often appears as a real-space renormalization group flow describing the dynamics of geometrical couplings. Typical examples are afforded by non-linear $\sigma $-model theory \cite{9} or by the averaging of cosmological spacetimes \cite{4,5}. In such a setting, in order to estimate the net effect of renormalization on the scaling of geometrical parameters, it is often desirable to establish monotonicity results for the various curvature functionals associated with the flow.  This is an extremely non-trivial task and the recent results by G. Perelman \cite{18} provide an important unexpected breakthrough of vast potential use in geometrical physics. Apparently inspired by the dilatonic action in string theory, Perelman has introduced \cite{18} the two functionals
\begin{equation*}
F[{g};f]\doteq \int_{\Sigma }({R}(\beta )+|\nabla f|^{2})e^{-f}d\mu _{{g}},
\label{F-funct}
\end{equation*}
and  
\begin{equation*}
{W}[{g};f_{\beta },{\tau }]\doteq \int_{\Sigma }\left[ {\tau} \left( \left| \nabla
f_{\beta }\right| ^{2}+{R}(\beta) \right) +f_{\beta }-3\right]\frac{e^{-f_{\beta }}}{(4\pi {\tau}
(\beta ))^{\frac{3}{2}}}\,\,d\mu _{{g}(\beta )},
\label{W-funct}
\end{equation*}
(see below for notation) depending on the geometry of the Riemannian manifold $(\Sigma ,g(\beta ))$ undergoing a Ricci flow evolution $\beta \rightarrow g(\beta )$,\,$\beta \in [0,T)$,  and on the choice of a probability measure $d\varpi (\beta )$ $=$
$(4\pi {\tau}(\beta ))^{-\frac{3}{2}}$\,$Vol[\Sigma _{\beta }]e^{-f_{\beta }}d\Pi _{\beta }$, ($d\Pi _{\beta }\doteq Vol^{-1}[\Sigma _{\beta }]d\mu _{{g}(\beta )}$ denoting the normalized Riemannian volume element), associated with a backward diffusion of the function $f$. The functional $F[{g};f]$ has, as already stressed, the structure of the dilatonic action, familiar in non-linear $\sigma $-model theory and in the statistical mechanics of extended objects. ${W}[{g};f_{\beta },{\tau }]$  is basically a scale-invariant generalization of $F[{g};f]$ associated with the introduction of the scale parameter $\tau (\beta )$, which controls the localization properties of the measure $d\varpi (\beta )$. The basic property of $F[{g};f]$ and ${W}[{g};f_{\beta },{\tau }]$ is their weakly monotonic character along Ricci flow trajectories, a fact that has been put to use  by Perelman \cite{18,19,20} in his work on the proof of Thurston geometrization conjecture. Such monotonicity properties and a few formal similarities with standard entropies in statistical mechanics  accounts for the attribute entropic. Further justifications  come from a closer look into the  structure of $F[{g};f]$ and ${W}[{g};f_{\beta },{\tau }]$. In this connection, the functional 
${W}[{g};f_{\beta },{\tau }]$ is particularly interesting since, as is easily checked, it contains a natural combination of the relative entropy $S[d\varpi (\beta )||d\Pi _{\beta }]$ associated with the pair of measures 
$(d\Pi  _{\beta },d\varpi (\beta ))$, of the corresponding entropy generating functional $I[d\varpi (\beta )||d\Pi _{\beta }]$, and of the $d\varpi (\beta )$-localized curvature average $<R(\beta )>_{d\varpi (\beta )}$$\doteq$$\int_{\Sigma }R(\beta )d\varpi (\beta )$. The particular form of their combination in ${W}[{g};f_{\beta },{\tau }]$ is strongly suggested by the theory of logarithmic Sobolev inequalities \cite{18,10}. Note that  $<R(\beta )>_{d\varpi (\beta )}$ enters as a defective parameter setting the size of scalar curvature over the region where $d\varpi (\beta )$ is localized.  It is interesting to remark that none of the constituents of 
${W}[{g};f_{\beta },{\tau }]$ has, by itself, any manifest monotonicity property along Ricci flow trajectories, and it is just their overall interaction in ${W}[{g};f_{\beta },{\tau }]$ that makes the shrinker entropy monotonic. Various authors \cite{8,13e,14,15} have exploited the strategy suggested by Perelman's construction and succeded in specializing or extending  ${W}[{g};f_{\beta },{\tau }]$ to other specific settings. However,  the natural question of the monotonicity of the constituent entropic functionals generating ${W}[{g};f_{\beta },{\tau }]$ does not seem to have received particular attention. Such an analysis is relevant to the physical applications of the Ricci flow and also for a deeper understanding of the properties of Perelman's shrinker entropy. In this paper we discuss such an issue in connection with the volume-normalized Ricci flow (such a choice being motivated by our long-standing interest in cosmological applications of the theory). Our main results are twofold. The analysis of the volume-normalized version of the Hamilton-Perelman flow easily shows
that by
a natural renormalization $\tau(\beta )\mapsto \widetilde{\tau }(\beta )$ of the scale parameter $\tau (\beta )$ 
one can make explicitly monotonic the curvature term $\tau (\beta )<R(\beta )>_{d\varpi (\beta )}$ appearing in ${W}[{g};f_{\beta },{\tau }]$. In this way, we can  connect the growth properties of scalar curvature to the localization behavior of $d\varpi (\beta )$.  This  strategy suggests also that, by deforming Perelman's backward diffusion $\{d\varpi(\beta ) \}_{\beta <T}$ $\mapsto $ $\{d\Omega (\beta ) \}_{\beta <T}$ by adding a suitable drift term, one may get also monotonicity for the corresponding renormalized relative entropy $S[d\Omega  (\beta )||d\Pi _{\beta }]$. Quite remarkably, the answer is in the affirmative and the resulting deformation is provided by a Fokker-Planck (backward) diffusion 
$\left\{d\Omega _{t}\right\}_{t\geq 0}$, $t\doteq \beta ^{*}-\beta$, with a drift term generated by the scalar curvature fluctuations. We show that, in a well-defined sense, this Fokker-Planck process is the natural diffusion along the (volume-normalized) Ricci flow. The analysis of $\left\{d\Omega _{t}\right\}_{t\geq 0}$ shows that we are dealing with a gradient flow generated by a (weakly) monotonic relative entropy $S[d\Omega _{t}||d\Pi _{t}]$. As an elementary consequence of such a monotonicity, we prove that if the Ricci curvature is positive then one gets exponential convergence of $\left\{d\Omega _{t}\right\}_{t\geq 0}$ to $d\Pi _{t}$.\\

\noindent  Our analysis relies on a remarkable parametrization of diffusion processes suggested by F. Otto \cite{16,17}  related to the use of the Wasserstein metric on measure spaces \cite{21a}, (see \cite{0} for an in depth and very informative presentation of the whole subject). The distance induced by such a metric provides a way of turning the space of probability measures on a Riemannian manifold into a geodesic space, and has recently drawn attention in attempts of extending the notion of Ricci curvature to general metric spaces \cite{13b,13c,21a,23}. The preliminary results presented in this paper point to the possibility, recently advocated also by J. Lott and C. Villani, that the use of the geometry of the space of probability measures and of the associated notions of optimal transport and Wasserstein metric may also play a significant role in Ricci flow theory. (\emph{Added in the arXive version v3}: important developments relating Wasserstein distance and Ricci flow theory have been recently considered also by R. J. McCann and P. Topping \cite{13d}, useful remarks in this connection are also discussed in the nice monography \cite{21ab} by P. Topping. The enlightening book by C. Villani \cite{22a} provides a most inspiring analysis of the deep interplay between optimal transportation theory and Riemannian geometry). 

\noindent \emph{Outline of the paper}. We star by recalling a few basic properties of the space of probability measures on Riemannian manifolds. In section 1.1 we discuss the geometry of such a space from the point of view advocated by F. Otto. In particular we analyze the case when a curve of probability measures covers a fiducial curve of Riemannian metrics. The natural framework for such a discussion is not a fixed probability space but rather a bundle of probability measure spaces over the space of Riemannian metrics, where each fiber is a probability space with a distinguished reference measure, (the normalized Riemannian volume element), metrized by the Wasserstein distance with cost function determined by the given Riemannian distance function. This is the situation occurring in the application of the formalism to Ricci flow theory. To our knowledge, such a general framework is not explicitly discussed in the existing literature on optimal transportation theory, and thus we pause a little bit for analyzing it in some detail. In such a setting we discuss explicitly the properties of gradient flows in the bundle space of probability measure. We conclude (section 1.2) this long overview of the probabilistic formalism by recalling the characterization of Wasserstein distance and its interplay with the relative entropy and the entropy production functionals. The connection with the formalism developed by Otto comes by when discussing the characterization of Wasserstein length of a curve of probability measures and its geodesic nature. Here we analyze in some depth the extension of the notion of Wasserstein length of a curve of probability measure to the case when the curve in question covers a reference curves of Riemannian metrics. Again, this case does not appear in the standard  literature and needs to be explicitly addressed. In particular, we emphasize that the characterization of Wasserstein geodesics as solution of a Hamilton--Jacobi equation cannot be trivially extended to this more general case. In our opinion this is a basic issue to be solved in order to apply optimal transportation theory to the Ricci flow. We comment on a possible approach to a strategy for a solution in the final part of the paper. In section 2 we discuss the Perelman coupling for the volume--normalized Ricci flow. Section 2.1 recalls a few properties of the shrinker entropy, some really well--known, and a few others not so easily spotted in the existing literature. In particular, by elaborating on a remark by A. Figalli, we explicitly show that the shrinker entropy  is an entropy balance functional basically generated by the time derivative of the relative entropy associated with Perelman backward heat equation. In section 2.2 we prove that there is a natural combination of a scale parameter and of the average curvature $\langle R \rangle_{d\varpi }$ which is weakly--monotonic along the Ricci flow. This curvature entropy gives rise to a useful $d\varpi $--averaged Harnack--type estimate. In section 3 we introduce the relation between Ricci flow and Fokker--Planck diffusion. This exploits Otto's parametrization of a probability measure by introducing a potential for scalar curvature fluctuations (section 3.1). Such a potential has a familiar counterpart in the Ricci flow theory for surfaces, and plays a distinguished role in our analysis. In particular we use it for estimating the Wasserstein length of the curve of normalized Riemannian volume elements along the Ricci flow. In section 3.2 we exploit the formalism so developed for discussing, under very general conditions, the evolution of an absolutely continuous curve of probability measure along the (backward) Ricci flow. In section 3.3 these results are used to prove that Fokker-Planck diffusion is the natural diffusion of a probability measure along the Ricci flow. In particular, the associated relative entropy is weakly--monotonic and the flow is  gradient--like. We also show that the flow of the associated Radon--Nikodym derivatives (with respect to the evolving Riemannian measure) is a true gradient flow with respect to Otto's inner product. We also emphasize the associated contraction properties in the corresponding (quadratic) Wasserstein distance, and their role in discussing the trend to equilibrium for the Fokker--Planck diffusion. In such a setting one naturally discovers that the associated relative entropy is displacement--convex. A property, this latter, which is strongly reminiscent of the characterization of Wasserstein geodesics in the space of probability measures. This point is discussed by showing that Fokker--Planck diffusion along the backward Ricci flow can be equivalently rewritten as a viscous Hamilton--Jacobi equation, where the viscosity parameter is related with the lower bound of the Ricci curvature. 
The paper concludes with an appendix stressing a few basic differences between Fokker--Planck diffusion and Perelman diffusion along the backward Ricci flow.

\section{Probability measures on Riemannian manifolds}
\label{Botto}
Throughout this paper, $\Sigma $ will denote a smooth three-dimensional manifold, which we assume to be closed and without boundary. We let $C^{\infty }(\Sigma, \mathbb{R})$ and $\doteq {C}^{\infty }(\Sigma,\otimes ^{2}_{+}\, T^{*}\Sigma)$ be the space of smooth functions and of smooth definite positive symmetric bilinear forms on $\Sigma $, respectively.  
$\mathcal{D}iff(\Sigma )$ is the group of smooth diffeomorphisms of $\Sigma $, and  $\mathcal{R}iem(\Sigma )$\,  $\doteq {C}^{\infty }(\Sigma,\otimes ^{2}_{+}\, T^{*}\Sigma)$ is the space of all smooth Riemannian metrics over $\Sigma$. 
The tangent space\,, $\mathcal{T}_{(\Sigma ,g)}\mathcal{R}iem(\Sigma )$, to $\mathcal{R}iem(\Sigma )$ at $(\Sigma,g )$ can be naturally identified with the space of symmetric bilinear forms ${C}^{\infty }(\Sigma,\otimes ^{2}\, T^{*}\Sigma)$ over $\Sigma $. It is endowed with the pre--Hilbertian  $L^{2}(\Sigma, g)$ inner product defined  on $(\Sigma,g)$ by 
\begin{equation}
(W,V)_{L^{2}(\Sigma )}\doteq \int_{\Sigma }g^{il}\,g^{km}\,W_{ik}\,V_{lm}d\mu _{g},
\label{L2prod}
\end{equation}
$W$, $V$ $\in {C}^{\infty }(\Sigma,\otimes ^{2}\, T^{*}\Sigma)$\ being (square-summable).
The hypothesis of smoothness has been made for simplicity. 
Results similar to those described below,  can be obtained for
finite H\" older or Sobolev differentiability.
\noindent In such a framework,
Let $d\mu _{g}$, $Vol\left[ \Sigma %
\right] _{g}$, and $d\Pi_{g}
\doteq Vol\left[ \Sigma \right] _{g}^{-1}d\mu _{g}$ respectively denote
the Riemannian density, the volume, and the corresponding
normalized measure on a Riemannian manifold $(\Sigma ,g)\in \mathcal{R}iem(\Sigma )$.
 In what follows, we will often refer (rather informally) to the bundle $%
\pi :Prob(\Sigma )\rightarrow \mathcal{R}iem(\Sigma )$ of all Borel probability measure on $\Sigma$, which are absolutely continuous with respect to normalized Riemannian volume element $d\Pi_{g}$. Each fiber $Prob(\Sigma,g)\doteq \pi ^{-1}(\Sigma,g )$ is endowed with the topology of weak 
convergence, and can be parametrized by the set of all Radon--Nikodym derivatives with respect to $d\Pi_{g} $, \emph{i.e.}, 
\begin{equation}
Prob(\Sigma,g )\doteq \left\{N\,d\Pi _{g}\,:\, N\in C_{b}(\Sigma ,\mathbb{R}^{+}),\,\int_{%
\Sigma }N\,d\Pi_{g} =1\right\} ,
\label{unouno}
\end{equation}
where $C_{b}(\Sigma ,\mathbb{R}^{+})$ is the space of positive bounded
measurable functions \cite{7}, (again, we often restrict our analysis to the smooth functions in $C_{b}(\Sigma ,\mathbb{R}^{+})$). To avoid notational prolixity, given a probability measure $d\varpi $ on $(\Sigma,g) $, we shall write for simplicity $d\varpi \in Prob(\Sigma, g )$ to actually mean $\left(\frac{d\varpi }{d\Pi_{g} }\right)\,d\Pi_{g}\in Prob(\Sigma,g )$. Moreover, for later convenience we shall restrict our
attention to smooth probability measure with finite $k$-th moments, $k\geq 1$, \emph{i.e.}, we assume, (but this is not strictly necessary as long as $(\Sigma ,g)$ is compact), that the Radon--Nikodym derivative $N$ satisfies $\int_{\Sigma }[d_{g}(x,y)]^{k}Nd\Pi_{g} $ $<\infty $ for some (and hence all) $x\in (\Sigma ,g)$, where $d_{g}(x,y)$ denotes the Riemannian geodesic distance in $(\Sigma,g) )$. Typically we set $k=2$. 
\subsection{Otto's parametrization}
\label{Ottop}
As suggested by F. Otto \cite{16}, (see also the remarkable paper \cite{17} on which this section is based and from which I extracted many observations), when discussing  probability diffusion semigroups on a Riemannian manifold $(\Sigma ,g)$ it can be profitable to consider each fiber $Prob(\Sigma,g)\in Prob(\Sigma )$ as an infinite dimensional manifold
locally modelled over 
the Hilbert space completion of the tangent space 
\begin{equation}
T_{N}Prob(\Sigma,g )\doteq \left\{ h\in C_{b}(\Sigma ,\mathbb{R}%
),\,\int_{\Sigma }h\,Nd\Pi_{g} =0\right\} ,
\end{equation}
with respect to the inner product defined, at the given Radon--Nikodym derivative $N$, by the Dirichlet form
\begin{equation}
\left\langle \varphi , \zeta \right\rangle _{(g,N)}\doteq
\int_{\Sigma }\left(g^{ik}\,\nabla _{k}\varphi \,\nabla _{i}\zeta\right) \,N\,d\Pi_{g}\;,
\label{inner}
\end{equation}
for any $\varphi$, $ \zeta $ $\in C_{0}^{\infty }(\Sigma, \mathbb{R})$. (Recently, this matter has been discussed from a geometric point of view in a series of papers by J. Lott and C. Villani, \cite{13ab,13ac,13abc}). 
Under such an identification, one can represent vectors in \ $T_{N}Prob(\Sigma,g )$ as the solutions 
of an elliptic problem naturally associated with the given probability
measure $N\,d\Pi_{g} $ according to 
\begin{equation}
(h,N)\in T_{N}Prob(\Sigma,g )\times Prob(\Sigma,g )\longmapsto \psi \in
C_{b}(\Sigma ,\mathbb{R})/\mathbb{R} ,
\end{equation}
where, for any given pair $(h,N)$, the function $\psi $ is formally determined on the given $(\Sigma ,g)$ by the elliptic PDE
\begin{equation}
-\nabla ^{i}\left( N\,\nabla _{i}\psi \right) =h,
\label{otto}
\end{equation}
under the equivalence relation  identifying any two such 
solutions differing by an additive constant.  In general, such a characterization is somewhat heuristic, at least in the sense that its validity must be
checked case by case, (a particularly clear and deep analysis of the whole topic is discussed in \cite{0}). As we shall see, it applies in our setting, and it provides a useful framework for discussing the entropic aspects of the volume-normalized Ricci flow.

 \noindent There is a further aspect about the geometry of $Prob(\Sigma )$ which it will be useful to have at our disposal. 
First, note that the tangent space to the bundle $Prob(\Sigma )$ at $(g,\,N\,d\Pi _{g})$ can be decomposed as
\begin{equation}
T_{(g,\,N\,d\Pi _{g})}\,Prob(\Sigma )=T_{g}\,\mathcal{R}iem(\Sigma )\oplus T_{N}\,Prob(\Sigma ,g)\;,
\end{equation}
a decomposition which, since $\mathcal{R}iem(\Sigma )$ is contractible, extend to the whole tangent bundle $T\,Prob(\Sigma )$.
Let $\Gamma :[0,1]\ni \lambda\mapsto g_{ab}(\lambda )$ be a smooth curve of metrics in $\mathcal{R}iem(\Sigma )$, with $Vol\,(\Sigma ,g(\lambda ))$ $=$ $Vol\,(\Sigma ,g(\lambda=0 ))$, $\forall \lambda \in [0,1]$. By means of the corresponding one--parameter family of normalized volume elements $d\Pi_{\lambda }
\doteq Vol\left[ \Sigma \right] _{g(\lambda )}^{-1}d\mu _{g(\lambda )}$\, $\in Prob(\Sigma ,g(\lambda ))$, the curve  $\Gamma $ naturally lifts to a corresponding curve in the bundle $Prob(\Sigma )$,
\begin{eqnarray}
[0,1]&\longrightarrow & Prob(\Sigma )\\
\lambda &\longmapsto  & \left(g_{ab}(\lambda ),d\Pi _{\lambda}\right).\notag
\end{eqnarray}
The tangent vector $\in T_{(g(\lambda ),d\Pi _{\lambda })}\,Prob(\Sigma )$\ to 
such a curve at any given value of the parameter $\lambda $, say $\lambda =s$, can be readily characterized, in analogy with (\ref{otto}), if we parametrize  $\frac{%
\partial }{\partial \lambda }d\Pi _{\lambda }|_{\lambda =s}$ in terms of a potential $\Theta 
_{s}$ obtained as the solution of the elliptic equation 
\begin{equation}
g^{ik}(s)\nabla _{i}\left( d\Pi _{s}\,\nabla _{k}\Theta  _{s}\right) =-\left. 
\frac{\partial }{\partial \lambda }d\Pi _{\lambda }\right| _{\lambda =s}=-d\Pi _{s }\left. g^{ik}(\lambda )\frac{\partial }{\partial \lambda }\,g_{ik}(\lambda )\right|_{\lambda =s}\;,
\label{elliptictheta}
\end{equation}
where $\nabla $ denotes the covariant derivative with respect to the metric $g_{ab}(\lambda )$. Since $\frac{\partial }{\partial \lambda }\,g_{ik}(\lambda )$ $\in T_{g(\lambda )}\mathcal{R}iem(\Sigma )$, the equation (\ref{elliptictheta})  is formulated in the tangent space $T_{(g(\lambda ),d\Pi _{\lambda })}\,Prob(\Sigma )$, (this is the reason why we have  expressed it in terms of the measure density $d\Pi _{\lambda }$ rather than of the corresponding Radon--Nikodym derivative).  Note that $d\Pi _{\lambda }$ is covariantly constant over the corresponding $(\Sigma ,g(\lambda ))$, 
\begin{equation}
\nabla _{i}d\Pi _{\lambda  }=d\Pi _{\lambda  }\,g^{ab}(\lambda  )\nabla
_{i}g_{ab}(\lambda  )=0,
\label{covconst}
\end{equation}
(this is equivalent to the familiar formula $\partial _{i}\ln \sqrt{g(\lambda  )%
}=\delta _{a}^{c}\Gamma _{ic}^{a}(\lambda  )$, where $\Gamma _{ic}^{a}(\lambda  )$
are the Christoffel symbols associated with $g_{ab}(\lambda  )$), thus
 we can rewrite (\ref{elliptictheta}) as $\triangle\,\Theta 
_{s}=-\,g^{ik}(\lambda )\frac{\partial }{\partial \lambda }\,g_{ik}(\lambda )|_{s}$, where $\triangle \doteq g^{ik}(\lambda )\nabla_{i} \nabla_{k} $ denotes the Laplace--Beltrami operator on $(\Sigma ,g(\lambda ))$.

\noindent It must be stressed that the family of potentials $\Theta _{\lambda }$ depends on the chosen  curve of metrics $\lambda \mapsto g_{ab}(\lambda )$, and not only from the associated Riemannian measures $\lambda \mapsto d\Pi _{\lambda }$. Actually, the dependence from $d\Pi _{\lambda }$ can be easily traded for the action of the diffeomorphisms group $\mathcal{D}iff(\Sigma )$.  This follows by observing that since we have normalized  the Riemannian volume elements $d\Pi _{\lambda }$, we can apply  Moser's theorem 
(see \emph{e.g.} \cite{3}, 7.2.3) according to which, on a compact
manifold $\Sigma $ admitting two volume forms $d\mu $ and $d\nu $ with $%
\int_{\Sigma }d\mu =\int_{\Sigma }d\nu $, there exists a diffeomorphism $%
\phi :\Sigma \longrightarrow \Sigma $ such that $\phi ^{\ast }d\mu =d\nu $.
In our case,  this implies that there exists a $\lambda $-dependent diffeomorphism 
\begin{gather}
\phi _{\lambda }:(\Sigma ,g_{ab}(\lambda ))\longrightarrow (\Sigma ,g_{ab}(\lambda =0)) \label{diffphi} \\
y^{k}\longmapsto x^{i}=\phi _{\lambda }^{i}(y^{k},\lambda )  \notag
\end{gather}
such that 
\begin{equation}
d\Pi _{\lambda }(y^{k})=J(\phi _{\lambda })\,\,d\Pi _{\lambda =0}(x^{i}(y^{k}))
\label{moserJ}
\end{equation}
where 
\begin{equation}
J(\phi _{\lambda })\doteq \left| \frac{\partial x^{i}(y^{k},\lambda )}{\partial y^{h}}%
\right| ,
\end{equation}
is the Jacobian of $\phi _{\lambda }$. In terms of the diffeomorphism $\phi _{\lambda }$, and of the metric $g_{ik}(\lambda )$, we can rewrite (\ref{elliptictheta}) as 
\begin{equation}
g^{ik}(s)\nabla _{i}\,\nabla _{k}\Theta  _{s}\,=-\left. 
\frac{\partial }{\partial \lambda }\,\ln\,J(\phi  _{\lambda })\,\right| _{\lambda =s}\;.
\label{ellipticln}
\end{equation}
Formally, given a solution $\Phi _{s}$ of (\ref{elliptictheta}) and any ($\lambda $-independent) smooth function with compact support $\zeta \in C_{0}^{\infty }(\Sigma, \mathbb{R})$, we have
\begin{equation}
\left.\frac{d }{d \lambda }\int_{\Sigma_{\lambda } }\zeta(x) d\Pi _{\lambda }\right|_{s}=\int_{\Sigma_{s}
}g^{ik}(\lambda )\,\nabla _{i}\zeta(x) \,\nabla _{k}\Theta  _{s}\,d\Pi _{s}= \left\langle \zeta ,\Theta  _{s}\right\rangle_{s}\;,
\label{th} 
\end{equation} 
where $\langle\ldots\,,\ldots\rangle_{s}$ is a shorthand notation for $\langle\ldots\,,\ldots\rangle_{(g(s),N=1)}$. 
According to (\ref{inner}) the 
relation (\ref{th}) identifies $(g_{ik}(\lambda ),d\Pi _{\lambda })\longmapsto (\frac{\partial }{\partial \lambda }g_{ik}(\lambda ),\,\Theta  _{\lambda }) $
as the tangent vector to the curve $\lambda \rightarrow (g_{ik}(\lambda ),d\Pi _{\lambda })$. The  family of function $\Theta  _{\lambda }$ play  also a fundamental role in characterizing gradient flows in the bundle $Prob(\Sigma )$ when, as in our case, the inner product $\langle\ldots\,,\ldots\rangle_{s}$ varies. To discuss this point, let us consider an absolutely continuous curve of probability measures $[0,1]\ni \lambda \mapsto d\varpi _{\lambda }\in Prob(\Sigma )$, covering the fiducial curve $\lambda \mapsto d\Pi  _{\lambda }$, (\emph{i.e.}, $d\varpi _{\lambda }$ is absolutely continuous with respect to $d\Pi _{\lambda }$, $\forall \lambda \in [0,1]$). Let us denote by $\Psi  _{\lambda }$ the tangent vector to $\lambda \mapsto d\varpi _{\lambda }$ defined, in analogy with (\ref{elliptictheta}), by the elliptic PDE
\begin{equation}
g^{ik}(s)\nabla _{i}\left( d\varpi   _{s}\,\nabla _{k}\Psi  _{s}\right) =-\left. 
\frac{\partial }{\partial \lambda }d\varpi   _{\lambda }\right| _{\lambda =s}\;.
\label{ellipticpsi1}
\end{equation}
In terms of  the associated Radon--Nikodym derivatives $\frac{d\varpi _{\lambda }}{d\Pi _{\lambda }}$, connecting the fiducial curve of reference measures $\lambda \mapsto d\Pi _{\lambda }$ to the curve $\lambda \mapsto d\varpi  _{\lambda }$, we have, (from the identity $\frac{\partial }{\partial \lambda }d\varpi   _{\lambda }=\frac{\partial }{\partial \lambda }(\frac{d\varpi _{\lambda }}{d\Pi _{\lambda }}\,d\Pi _{\lambda })$),
\begin{equation}
\left(\frac{\partial }{\partial \lambda }+\nabla ^{i}\Theta  _{\lambda }\,\nabla _{i} \right)\left.\frac{d\varpi _{\lambda }}{d\Pi _{\lambda }}\right| _{\lambda =s}=
-\nabla ^{i}\left[\frac{d\varpi  _{s}}{d\Pi _{s}}\,\nabla _{i}(\Psi _{s}-\Theta  _{s})\right] \;,
\label{convDer}
\end{equation}
which is easily seen to be interpretable in the sense
\begin{eqnarray}
&&\frac{d}{d\lambda }\left.\int_{\Sigma }\,\zeta(x)\,\,\left(\frac{d\varpi _{\lambda }}{d\Pi _{\lambda }}\right) \,d\Pi_{\lambda }\right| _{\lambda =s} = \label{RDtangent3}\\
\notag\\
&=&\int_{\Sigma }
\nabla _{i}\zeta(x)\, \nabla ^{i}(\Psi _{s}-\Theta  _{s})\,\left(\frac{d\varpi  _{s}}{d\Pi _{s}}\right) \,d\Pi _{s}\;,
\notag
\end{eqnarray}
$\forall \zeta(x) \in C_{0}^{\infty }(\Sigma )$. Namely, $(\Psi _{\lambda }-\Theta _{\lambda })$ is the tangent vector, at $\lambda =s$, to the curve of Radon--Nikodym derivatives  $[0,1]\ni \lambda \mapsto \left(\frac{d\varpi _{\lambda }}{d\Pi _{\lambda }}\right)\in C_{b}(\Sigma ,\mathbb{R}^{+}) $. 
If we evaluate the inner product between  $(\Psi _{\lambda }-\Theta _{\lambda })$ and a generic tangent vector $\xi\in T_{\frac{d\varpi _{\lambda }}{d\Pi _{\lambda }}}Prob(\Sigma ,g(\lambda ))$, then we get
\begin{eqnarray}
&&\left\langle (\Psi _{\lambda }-\Theta _{\lambda }),\,\xi  \right\rangle_{\frac{d\varpi _{\lambda }}{d\Pi _{\lambda }}}=\label{OL2}\\
&=&\int_{\Sigma _{\lambda }}\left[g^{ik}(\lambda )\,\nabla _{i}(\Psi _{\lambda }-\Theta _{\lambda })\,\nabla _{k}\xi\right] \;\left(\frac{d\varpi _{\lambda }}{d\Pi _{\lambda }} \right)d\Pi _{\lambda }\notag\\
&=&-\int_{\Sigma _{\lambda }}\xi\,\, \nabla _{i}\left[\frac{d\varpi _{\lambda }}{d\Pi _{\lambda }}\,g^{ik}(\lambda )\,\nabla _{k}(\Psi _{\lambda }-\Theta _{\lambda }) \right]d\Pi _{\lambda }\notag\\
&=&\int_{\Sigma _{\lambda }}\xi\,\left[\left(\frac{\partial }{\partial \lambda }+g^{ik}(\lambda )\,\nabla _{i}\Theta  _{\lambda }\,\nabla _{k} \right)\frac{d\varpi _{\lambda }}{d\Pi _{\lambda }}\right]\;d\Pi _{\lambda }\notag\\
&=& \left\langle \left(\frac{\partial }{\partial \lambda }+\,g^{ik}(\lambda )\,\nabla _{i}\Theta  _{\lambda }\,\nabla _{k} \right)\frac{d\varpi _{\lambda }}{d\Pi _{\lambda }} ,\,\xi  \right\rangle_{L^{2}({d\Pi _{\lambda }})} \;\notag,
\end{eqnarray}
where $\langle\,,\,\rangle_{L^{2}({d\Pi _{\lambda }})}$ denotes the standard ${L^{2}({d\Pi _{\lambda }})}$ inner product on $(\Sigma , g(\lambda ))$. Similarly, if we denote by $Grad\,\mathcal{F}$ the gradient of a smooth functional $\mathcal{F}:Prob(\Sigma )\longrightarrow \mathbb{R}$, with respect to the inner product  $\langle\ldots,\ldots\rangle_{\frac{d\varpi _{\lambda }}{d\Pi _{\lambda }}}$, then we compute
\begin{eqnarray}
&&\left\langle Grad\, \mathcal{F},\,\xi  \right\rangle_{\frac{d\varpi _{\lambda }}{d\Pi _{\lambda }}}=\label{OL2a}\\
&=&\int_{\Sigma _{\lambda }}\,g^{ik}(\lambda )\nabla _{i}\,Grad\, \mathcal{F} \,\,\nabla _{k}\xi \;\left(\frac{d\varpi _{\lambda }}{d\Pi _{\lambda }} \right)d\Pi _{\lambda }\notag\\
&=&-\int_{\Sigma _{\lambda }}\xi\,\, \nabla _{i}\left[\frac{d\varpi _{\lambda }}{d\Pi _{\lambda }}\,g^{ik}(\lambda )\,\nabla _{k}Grad\, \mathcal{F}
 \right]d\Pi _{\lambda }\notag\\
&=&-\left\langle \nabla _{i}\left[\frac{d\varpi _{\lambda }}{d\Pi _{\lambda }}\,g^{ik}(\lambda )\,\nabla _{k}Grad\, \mathcal{F}
 \right] ,\,\xi  \right\rangle_{L^{2}({d\Pi _{\lambda }})}\;.\notag
\end{eqnarray}
Note that both in (\ref{OL2}) and (\ref{OL2a}) we have explicitly inserted the $\lambda $--dependent metric $g^{ik}(\lambda )$ in order to make it clear that the inner product $\langle\ldots,\ldots\rangle_{\frac{d\varpi _{\lambda }}{d\Pi _{\lambda }}}$ depends from the curve of  Riemannian metrics $\lambda \mapsto g_{ab}(\lambda )$.

\noindent When $\mathcal{F}$ is identified with  the relative entropy functional \cite{7}, 
\begin{equation}
S\left[ d\varpi \parallel d\Pi_{g}  \right] \doteq \left\{ 
\begin{tabular}{lll}
$\int_{\Sigma }\frac{d\varpi }{d\Pi_{g}  }\ln \frac{d\varpi }{d\Pi_{g}  }d\Pi_{g}  $ & if & 
$d\varpi <<d\Pi_{g}  $ \\ 
&  &  \\ 
$\infty $ &  & otherwise,
\end{tabular}
\right. 
\label{entfunct}
\end{equation}
where $d\varpi <<d\Pi_{g}  $ stands for absolute continuity,  the relations (\ref{OL2}) and (\ref{OL2a}) allow to characterize a class of flows $\lambda \mapsto \frac{d\varpi _{\lambda }}{d\Pi _{\lambda }}$ which will be important in what follows. We start computing $Grad\,S[d\varpi\parallel d\Pi_{g}]$ in such a case. To this end,  let us consider a linearization of $\frac{d\varpi}{d\Pi_{g}}$ in the direction of the generic vector $\xi \in T_{\frac{d\varpi}{d\Pi_{g}}}\,Prob(\Sigma ,g)$,  \emph{i.e.},
\begin{equation}
\frac{d\varpi}{d\Pi_{g}}(\epsilon )\doteq \frac{d\varpi}{d\Pi_{g}}+\epsilon \,\,\xi\;, 
\end{equation}
with $\xi $ parametrized \emph{a' la Otto},
\begin{equation}
\xi =-\nabla ^{i}\left[\frac{d\varpi}{d\Pi_{g}}\,\nabla _{i}\Upsilon  \right]\;.
\end{equation}
The corresponding linearization of $S[d\varpi\parallel d\Pi_{g}]$ in the direction $\xi $ provides
\begin{eqnarray}
&& D\,S[d\varpi\parallel d\Pi_{g}]\circ\xi \doteq  \left.\frac{d}{d\,\epsilon }\,S_{\epsilon }[d\varpi\parallel d\Pi_{g}]\right|_{\epsilon =0} =\\
&=& -\int_{\Sigma }\left(1+\ln\,\frac{d\varpi}{d\Pi_{g}} \right)\nabla ^{i}\left[\frac{d\varpi}{d\Pi_{g}}\,\nabla _{i}\Upsilon  \right]d\Pi_{g}\notag\\
&=& \int_{\Sigma }\nabla ^{i}\,\ln\,\frac{d\varpi}{d\Pi_{g}}\,\nabla _{i}\Upsilon\,d\Pi_{g}\notag\\
&=&\left\langle \ln\,\frac{d\varpi}{d\Pi_{g}},\,\,\Upsilon \right\rangle_{\frac{d\varpi}{d\Pi_{g}}}\notag\;,
\end{eqnarray}
which implies that the gradient, with respect to $\langle\ldots,\,\ldots\rangle_{{d\varpi}\backslash {d\Pi_{g}}}$, is given by 
\begin{equation}
Grad\, S[d\varpi\parallel d\Pi_{g}]=\,\ln\,\frac{d\varpi}{d\Pi_{g}}\;.
\label{gradS1}
\end{equation}

\noindent With these preliminary remarks along the way, given the flow of reference Riemannian metrics $\lambda \mapsto g_{ab} ({\lambda })$, let us consider  a curve  $[0,1]\ni \lambda \mapsto \frac{d\varpi _{\lambda }}{d\Pi _{\lambda }}\in C^{\infty}(\Sigma ,\mathbb{R})\cap C_{b}(\Sigma ,\mathbb{R}^{+})$ whose tangent vector  $(\Psi _{\lambda }-\Theta _{\lambda })$ is such that

\begin{equation}
\left\langle (\Psi _{\lambda }-\Theta _{\lambda }),\,\xi  \right\rangle_{\frac{d\varpi _{\lambda }}{d\Pi _{\lambda }}}+
\left\langle Grad\, S[d\varpi\parallel d\Pi_{g}],\,\xi  \right\rangle_{\frac{d\varpi _{\lambda }}{d\Pi _{\lambda }}}=0\;,
\label{gradcond}
\end{equation}

\noindent $\forall \xi \in T_{\frac{d\varpi _{\lambda }}{d\Pi _{\lambda }}}Prob(\Sigma ,g(\lambda ))$. Note that, according  to (\ref{OL2}) and (\ref{OL2a}), such a condition is equivalent, in the ${L^{2}(\Sigma ,{d\Pi _{\lambda }})}$  sense, to the PDE

\begin{equation}
\left(\frac{\partial }{\partial \lambda }+\,g^{ik}(\lambda )\,\nabla _{i}\Theta  _{\lambda }\,\nabla _{k} \right)\frac{d\varpi _{\lambda }}{d\Pi _{\lambda }}=
\nabla _{i}\left[\frac{d\varpi _{\lambda }}{d\Pi _{\lambda }}\,g^{ik}(\lambda )\,\nabla _{k}Grad\, S
 \right]\;,
 \label{GFG}
\end{equation}

\noindent where $Grad\, S$ is a shorthand notation for $Grad\, S[d\varpi\parallel d\Pi_{g}]$. If we insert in this latter condition the expression (\ref{gradS1}) for $Grad\,S$, we get 
the Fokker--Planck equation
\begin{equation}
\left(\frac{\partial }{\partial \lambda }+\nabla ^{i}\Theta  _{\lambda }\,\nabla _{i} \right)\frac{d\varpi _{\lambda }}{d\Pi _{\lambda }}=
\triangle \left(\frac{d\varpi _{\lambda }}{d\Pi _{\lambda }}
 \right)\;.
 \label{GFS}
\end{equation}

\noindent An elementary computation shows that
\begin{eqnarray}
\frac{d}{d\lambda }\,S[d\varpi_{\lambda }\parallel d\Pi_{\lambda }]=\left\langle (\Psi _{\lambda }-\Theta _{\lambda }),
\,S[d\varpi_{\lambda }\parallel d\Pi_{\lambda }] \right\rangle_{\frac{d\varpi _{\lambda }}{d\Pi _{\lambda }}}\;.
\end{eqnarray}
Thus, if we set $\xi =S[d\varpi_{\lambda }\parallel d\Pi_{\lambda }]$  in (\ref{gradcond}), we get

\begin{equation}
\frac{d}{d\lambda }\,S[d\varpi_{\lambda }\parallel d\Pi_{\lambda }]=-\left\langle Grad\, S,\,Grad\, S \right\rangle_{\frac{d\varpi _{\lambda }}{d\Pi _{\lambda }}}\;,
\label{gradcond2}
\end{equation}

\noindent which implies that (\ref{gradcond}) is the condition for the curve $[0,1]\ni \lambda \mapsto \frac{d\varpi _{\lambda }}{d\Pi _{\lambda }}\in C^{\infty}(\Sigma ,\mathbb{R})\cap C_{b}(\Sigma ,\mathbb{R}^{+})$ to be the gradient flow of the relative entropy functional $S[d\varpi_{\lambda }\parallel d\Pi_{\lambda }]$, with respect to the $\lambda $--varying inner product $\langle,\,\,\rangle_{{d\varpi}\backslash {d\Pi_{g}}}$. Note that, in the $L^{2}(\Sigma ,d\Pi _{g})$--sense, such a gradient flow condition implies that the Radon--Nikodym derivatives $\frac{d\varpi_{\lambda }}{d\Pi_{\lambda }}$ evolves  according to the Fokker--Planck diffusion (\ref{GFS}).

\noindent The relation (\ref{gradcond}) is an elementary but important property of the evolution of the relative entropy functional $S[d\varpi_{\lambda }\parallel d\Pi_{\lambda }]$ along a fiducial curve of Riemannian metrics $\lambda \mapsto g_{ab}(\lambda )$. It is ultimately related to Moser's theorem. To disclose the rationale underlying this latter remark, let us recall the well--known fact that, for a given fixed metric $g$, the gradient flow $\left(d\varpi_{\lambda },d\Pi_{g}\right)_{\lambda \geq 0}$ of $S[d\varpi_{\lambda }\parallel d\Pi_{g}]$ is the standard heat flow on $(\Sigma ,g)$, \cite{21b,22,22a}.  Along the same vein, let us observe that, along  a fiducial curve of Riemannian metrics $\lambda \mapsto g_{ab}(\lambda )$, we can write
\begin{eqnarray}
&& S[d\varpi_{\lambda }\parallel d\Pi_{\lambda }]=\int_{\Sigma }\frac{d\varpi_{\lambda } }{d\Pi _{\lambda }}\,
\ln\,\frac{d\varpi_{\lambda } }{d\Pi _{\lambda }}\,d\Pi_{\lambda }\\
\notag\\
&=&\int_{\Sigma }\frac{d\varpi_{\lambda } }{d\Pi _{\lambda =0}}\,\ln\,\frac{d\varpi_{\lambda } }{d\Pi _{\lambda =0}}\,d\Pi_{\lambda =0}-
\int_{\Sigma }\frac{d\varpi_{\lambda } }{d\Pi _{\lambda =0}}\,\ln\,\frac{d\Pi _{\lambda } }{d\Pi _{\lambda =0}}\,d\Pi_{\lambda =0}\notag\\
\notag\\
&=&S[d\varpi_{\lambda }\parallel d\Pi_{\lambda=0 }]-\int_{\Sigma }\frac{d\varpi_{\lambda } }{d\Pi _{\lambda =0}}\,\ln\,(J(\phi _{\lambda }))\,d\Pi_{\lambda =0}\notag\;,
\end{eqnarray}
where we have exploited Moser's theorem in the form (\ref{moserJ}). Thus, the relative entropy functional  $S[d\varpi_{\lambda }\parallel d\Pi_{\lambda }]$ of a diffusion process $\left(d\varpi_{\lambda }\right)_{\lambda \geq =0}$, with respect to a $\lambda $--varying reference measure $\left(d\Pi _{\lambda }\right)_{\lambda \geq =0}$, is the sum of the relative entropy of $\left(d\varpi_{\lambda }\right)_{\lambda \geq =0}$ with respect to the fixed Riemannian volume element $d\Pi _{\lambda =0},$ plus a forcing potential term provided by $\ln\,(J(\phi _{\lambda }))$. It is well--known that the gradient flows of relative entropies with forcing potentials typically yield for a Fokker-Planck diffusion \cite{22,22a}. Obviously, this heuristic explanation must be taken with care, because of the $\lambda $--dependence in the potential term  $\ln\,(J(\phi _{\lambda }))$. Nonetheless, it provides a natural framework for understanding the subtle interplay between the diffusion of probabilty measure on Riemannian manifolds  evolving along a  geometric flow $\lambda \mapsto (\Sigma ,g(\lambda ))$. Indeed, 
a central theme of this paper is that  Fokker--Planck dynamics has  remarkable geometric properties exactly when the fiducial curve of reference measures  $\lambda \mapsto d\Pi  _{\lambda }$ is generated by the (backward) Ricci flow.

\medskip 

\subsection{Relative entropy and the Wasserstein distance}

As we have seen above, the relative entropy functional  $S[d\varpi_{\lambda }\parallel d\Pi_{\lambda }]$, defined by (\ref{entfunct}), has a distinguished role in disclosing the interplay between Riemannian geometric flows and diffusion processes $(d\varpi_{\lambda })_{\lambda \geq 0}$. In particular, the approach to relative equilibrium $(d\varpi_{\lambda })_{\lambda \geq 0} \Rightarrow d\Pi_{g} $ is often  controlled by a logarithmic Sobolev inequality (LSI) \cite{7,10,17,21b}, that can be conveniently expressed \cite{7,17}, in terms of $S[d\varpi_{\lambda }\parallel d\Pi_{\lambda }]$, as  
\begin{equation}
LSI(\rho;B ) \doteq
\frac{1}{2\rho }(I[d\varpi \parallel d\Pi_{g}  ]+B)- S[d\varpi \parallel d\Pi_{g} ]\geq 0,
\label{logsob}
\end{equation}
where $\rho >0$ and  $B\geq 0$ are constants depending on the underlying geometry of $(\Sigma ,g)$,  and
\begin{equation}
I[d\varpi \parallel d\Pi_{g} ]\doteq \int_{\Sigma }\left|
\nabla \ln \frac{d\varpi}{d\Pi_{g} }\right| ^{2}d\varpi 
\end{equation}
is the entropy production functional (or Fisher information \cite{7}). In general we can set $B=0$ in (\ref{logsob}), however, as will become apparent after equation (\ref{defectiveW}) below, the defective form (\ref{logsob}) has some notational advantages in our setting. Recall also that if we assume the Bakry-Emery criterion \cite{2} $Ric(g)-Hess(\ln \frac{d\varpi }{d\Pi_{g}})$ $\geq$
$\rho\, {g}$, where $Hess(\circ )$  denotes the  Hessian on $(\Sigma ,g)$, then $LSI(\rho;B )$
holds on compact Riemannian manifolds without boundary  \cite{3a,4a,4b,7,21,14}. 

\noindent In this section we describe some of the  basic properties of $S[d\varpi_{\lambda }\parallel d\Pi_{\lambda }]$,  \cite{0}, \cite{17} that we shall need later on. Let us start by recalling that Jensen's inequality implies that  $%
S[d\varpi \parallel d\Pi_{g}  ]\in $ $[0,+\infty ]$, (this can also be checked directly by noticing that $S[d\varpi \parallel d\Pi_{g}  ]$ can be rewritten \cite{17}
as the integral of the non-negative function $
\frac{d\varpi }{d\Pi_{g}  }\left( \ln \frac{d\varpi }{d\Pi_{g}  }-1\right) +1$). Moreover, as a function of the
probability measures $d\varpi $ and $d\Pi_{g} $, the functional $S[d\varpi
\parallel d\Pi_{g}  ]$ is convex and lower semicontinuous in the 
weak topology on $Prob(\Sigma,g) $, and $S[d\varpi \parallel
d\Pi_{g}  ]=0$ iff $d\varpi =$ $d\Pi_{g}  $. It can be characterized \cite{7} by the
variational formula 
\begin{gather}
S[d\varpi \parallel d\Pi_{g}  ]= \\
=\sup_{f}\left\{ \int_{\Sigma }f\;d\varpi -\ln \int_{\Sigma }\exp
[f\;]\;d\Pi_{g}  \;:\,f\in C_{b}(\Sigma ;\mathbb{R})\right\} .  \notag
\end{gather}
Roughly speaking $S[d\varpi
\parallel d\Pi_{g}  ]$ provides the rate functional for the large deviation
principle \cite{7} controlling how deviant is the distribution of $d\varpi $ with
respect to the reference measure $d\Pi_{g}  $. In particular
one has \cite{7} (Pinsker's inequality) 
\begin{equation}
S[d\varpi \parallel d\Pi_{g} ]\geq \frac{1}{2}\left\| d\varpi -d\Pi_{g}  \right\|
_{var}^{2},
\label{pinsker}
\end{equation}
where we have introduced the total variation norm on $Prob(\Sigma,g )$ defined
by 
\begin{gather}
\left\| d\varpi -d\Pi_{g}  \right\| _{var}\doteq  \label{variation} \\
\doteq \sup_{\left\| f\right\| _{b}\leq 1}\left\{ \left| \int_{\Sigma
}fd\varpi -\int_{\Sigma }fd\Pi_{g}  \right| \,:\,f\in C_{b}(\Sigma ;\mathbb{R}%
)\,\right\} ,  \notag
\end{gather}
with $\left\| f\right\| _{b}\leq 1$ the uniform norm on $C_{b}(\Sigma ;%
\mathbb{R})$. This is a particular (and elementary) case of  transportation inequalities involving $S[d\varpi \parallel d\Pi_{g} ]$ and the 
notion of Wasserstein distance between probability measures \cite{0}, \cite{22}. Let us recall that for $d\varpi_{1}$, $d\varpi_{2}$ $\in $ $Prob(\Sigma,g)$,  we define the 
Wasserstein distance of order $s$ between $d\varpi_{1}$ and $d\varpi_{2}$ as
\begin{equation}
D_{s}^{W}(d\varpi_{1},\,d\varpi_{2})\doteq \inf_{\pi\in H(d\varpi_{1},\,d\varpi_{2}) }\left( \iint_{\Sigma \times
\Sigma }d(x,y)^{s}\pi (d\varpi_{1}\,\,d\varpi_{2})\right) ^{\frac{1}{s}}
\label{wassert}
\end{equation}
where $H(d\varpi_{1},\,d\varpi_{2})$ $\subset$  $Prob((\Sigma,g) \times (\Sigma,g) )$ denotes the set of probability
measures on $\Sigma \times \Sigma $ with marginals 
$d\varpi_{1}$ and $d\varpi_{2}$, \emph{i.e.}, such that $\pi (U \times \Sigma
)=d\varpi_{1} (U )$ and $\pi (\Sigma \times U )=d\varpi _{2} (U )$ for
any measurable set $U \subset \Sigma $; ($H(d\varpi_{1} ,d\varpi _{2} )$ is often called the set of couplings between $d\varpi_{1} $ and $d\varpi _{2} $). Note that by Kantorovich-Rubinstein duality, we have that $D_{1}^{W}(d\varpi_{1} ,d\varpi _{2} )$ $=$ $\left\| d\varpi_{1} -d\varpi _{2}  \right\| _{var}$.
Intuitively, $D_{s}^{W}(d\varpi_{1} ,d\varpi _{2} )$ represents, as we consider all possible couplings between the measures $d\varpi_{1}$ and $d\varpi _{2} $, the minimal cost needed to transport $d\varpi_{1}$ into $d\varpi _{2} $ provided that the cost to transport the point $x$ into the point $y$ is given by $d(x,y)^{s}$.
The distance
$D_{s}^{W}(d\varpi_{1} ,d\varpi _{2} )$ metrizes $Prob(\Sigma,g )$ turning it into a geodesic space. 

\noindent The pair 
$(Prob(\Sigma, g ),D_{2}^{W})$ has recently drawn attention \cite{13b,13c,15a, 21a,23} as an appropriate setting for extending the notion of Ricci curvature to general metric spaces. In this connection, a particularly elegant approach has been introduced in \cite{23}, by relating the ($K$--)convexity of the entropy functional $S[d\varpi \parallel d\Pi_{g} ]$ to Ricci curvature lower bounds. 
Explicitly, (\cite{23}, Th.1), if \ $\gamma :\,[0,1]\ni \lambda \rightarrow d\varpi _{\lambda }\in Prob(\Sigma,g )$ is a $(Prob(\Sigma, g ),D_{2}^{W})$-geodesic, (note that,
typically, $\gamma $ is not the linear interpolation of \ $d\varpi_{\lambda } $ and $%
d\Pi_{g} $), then a lower bound on  $Ric(g)$ is
equivalent to the $K$-convexity of $S[d\varpi_{\lambda } \parallel d\Pi_{g} ]$ along any
such geodesic $\gamma $, \emph{i.e.}
\begin{gather}
S[d\varpi_{\lambda }\parallel d\Pi_{g} ]\leq (1-\lambda  )S[d\varpi _{0}\parallel d\Pi_{g}
]+\lambda  S[d\varpi _{1}\parallel d\Pi_{g} ]- \label{KK} \\ 
-\frac{K}{2}\lambda  (1-\lambda  )[D_{2}^{W}(d\varpi _{0},d\varpi _{1})]^{2}  \notag
\end{gather}
iff $Ric(g)\geq K$, with $K\in \mathbb{R}$.
Such a result again points to transportation inequalities \cite{3b, 3c,13a,22}. 
\ For our purposes, it is sufficient
to recall the (simpler) bound (\cite{3c} case 5 of Th.1) which  always holds on compact Riemannian manifolds
\begin{equation}
D_{s}^{W}(d\varpi ,d\Pi_{g} )\leq 2^{\frac{1}{2s}}diam(\Sigma ,g)\,\,S[d\varpi
\parallel d\Pi_{g} ]^{\frac{1}{2s}},
\label{Pinsk2}
\end{equation}
$\forall\,d\varpi\in Prob(\Sigma, g )$, and
where $diam(\Sigma ,g)\doteq \sup \{d_{g}(x,y);\,x,y\in (\Sigma ,g)\}$
denotes the diameter of $(\Sigma ,g)$. It should be stressed that when, as in our case, one has a family of Riemannian manifolds $\lambda \mapsto (\Sigma ,g(\lambda ))$, such a Talagrand--like  inequality is effective as long as one has some uniform control on  $diam(\Sigma ,g(\lambda ))$.

\noindent A rather direct connection between Otto's description of $Prob(\Sigma ,g)$, discussed in the previous paragraph, and the (quadratic) Wasserstein distance $D_{2}^{W}(d\varpi_{1} ,d\varpi _{2})$ has been stressed by Otto and Villani \cite{17}, (see also \cite{0} for a more general setting), by relating $D_{2}^{W}(d\varpi_{1} ,d\varpi _{2} )$ to the geodesic distance associated with the inner product (\ref{inner}). This relation has been analyzed in detail by Lott \cite{13ac}, who has proved the following, (see \cite{13ac}, \ Prop. 3.3, and Prop. 4.24),

\begin{theorem} (Lott) Let $\varpi :\,[0,1]\ni \lambda \mapsto d\varpi _{\lambda }\in Prob(\Sigma ,g)$ be a smooth curve in $Prob(\Sigma ,g)$ with tangent vector $\Psi _{\lambda }$, defined for any fixed $\lambda=s$ by the elliptic PDE
\begin{equation}
g^{ik}\nabla _{i}\left( d\varpi   _{s}\,\nabla _{k}\Psi  _{s}\right) =-\left. 
\frac{\partial }{\partial \lambda }d\varpi   _{\lambda }\right| _{\lambda =s}\;,
\label{psilott}
\end{equation}
and such that $\nabla \Psi _{\lambda }\not=0$,\,$\forall \lambda \in [0,1]$, (\emph{i.e.}, the curve is immersed). Let $0=\lambda _{0}\leq \lambda _{1}\leq \ldots\leq \lambda _{J}=1$ be a partition of $[0,1]$, and let
\begin{equation}
L(\varpi )\doteq \sup_{J\in \mathbb{N}}\;\;\sup_{0=\lambda _{0}\leq \lambda _{1}\leq \ldots\leq \lambda _{J}=1}\;
\;\sum_{j=1}^{J}\;D_{2}^{W}(d\varpi _{\lambda _{j-1}},\,\,d\varpi _{\lambda _{j}})\;,
\label{walength}
\end{equation}
denote the length of the curve $\varpi $ in the Wasserstein space $\left(Prob(\Sigma ,g),\,D_{2}^{W} \right)$.
Then 
\begin{eqnarray}
 L(\varpi )&=&\int_{0}^{1}\left(\left\langle\,\Psi _{\lambda },\,\Psi _{\lambda }\,\right\rangle_{d\varpi  _{\lambda }}\right)^{\frac{1}{2}}\,\,d\lambda  \label{OttoWass}\\
&=& \int_{0}^{1}\left(\int_{\Sigma }\nabla ^{i}\Psi _{\lambda }\,\nabla _{i}\Psi _{\lambda }\,d\varpi  _{\lambda } \right)^{\frac{1}{2}}d\lambda\;.\notag
\end{eqnarray}
Moreover the curve $\varpi $ is a geodesic in $\left(Prob(\Sigma ,g),\,D_{2}^{W} \right)$ if its tangent vector $\Psi _{\lambda }$ satisfies the Hamilton-Jacobi equation
\begin{eqnarray}
\frac{\partial \Psi _{\lambda } }{\partial \lambda }\,+\,\frac{|\nabla \Psi _{\lambda }|^{2}}{2}=0\;,
\label{HamJac}
\end{eqnarray}
modulo the addition of a spatially--constant function to $\Psi _{\lambda }$.
\label{LottLength}
\end{theorem}

\medskip

\noindent It is important to discuss if the relation (\ref{OttoWass}) of Theorem \ref{LottLength} still holds if the metric $g_{ab}$, defining the reference measure and the diffusion operator (\ref{psilott}) in $Prob(\Sigma ,g)$, is replaced with a curve of (volume preserving) fiducial metrics $[0,1]\ni \lambda \mapsto g_{ab}(\lambda )$, with $g_{ab}(\lambda )$  uniformly bounded above and below for $0\leq \lambda \leq 1$. The potentially delicate issue concerns the characterization, along the given fiducial curve of metrics $\lambda \mapsto g_{ab}({\lambda })$,  of the Wasserstein distance (the  distance associated with the inner product (\ref{inner}) extends in an obvious manner along $\lambda \mapsto g_{ab}({\lambda })$). Indeed, the Wasserstein distance is usually defined in a fixed $Prob(\Sigma,g )$, (see (\ref{wassert})), with a cost function $d(\;,\;)^{2}$ provided by a fixed metric tensor $g_{ab}$. To characterize its extension to curves in $Prob(\Sigma )$, covering $\lambda \mapsto g_{ab}({\lambda })$, let us assume, in line with the boundedness hypotheses on $g_{ab}({\lambda })$, that there are constants $C>0$, (typically depending only on the dimension of $(\Sigma ,g(\lambda ))$), and $M>0$, (depending on a uniform bound on the geometry of $(\Sigma ,g(\lambda ))$),  such that
\begin{equation}
e^{-C\,M\,(\lambda _{j}-\lambda _{k})}\,d_{\lambda _{k}}(x,y)\leq d_{\lambda _{j}}(x,y)\leq e^{C\,M\,(\lambda _{j}-\lambda _{k})}\,d_{\lambda _{k}}(x,y)
\label{dbound} 
\end{equation}
for any points $x$ and $y$ in $\Sigma $, and any $\lambda _{k}$, $\lambda _{j}$ $\in [0,1]$, with $\lambda _{k}<\lambda _{j}$, (in the  Ricci flow case, (\ref{dbound}) holds whenever the Ricci curvature $Ric(g(\lambda ))$ of $(\Sigma ,g(\lambda ))$ is bounded by  $M$, $|Ric(g(\lambda ))|\leq M$, \cite{12}).  For any  $\lambda \in [0,\,1]$,  denote by 
\begin{equation}
\phi _{\lambda,\lambda _{k}}\,:\Sigma _{\lambda }\longrightarrow \Sigma _{\lambda _{k}}\;,
\end{equation}
the diffeomorphism defined by Moser's theorem, (see (\ref{diffphi})),  such that $d\Pi _{\lambda }=J(\phi _{\lambda,\lambda _{k}})\,d\Pi _{\lambda_{k} }$, where $J(\phi _{\lambda,\lambda _{k}})$ is the Jacobian of $\phi _{\lambda,\lambda _{k}}$. In line with the above hypotheses on the metric $g_{ab}(\lambda )$, we assume that there are  constants $C'$, $M'>0$, (again depending from a uniform bound on the geometry of $(\Sigma ,g(\lambda ))$),  such that
\begin{equation}
e^{-C'\,M'\,(\lambda -\bar \lambda)}\,J(\phi _{\bar \lambda,\lambda _{k}})\leq J(\phi _{\lambda,\lambda _{k}})\leq e^{C'\,M'\,(\lambda _{j}-\bar \lambda)}\,J(\phi _{\bar\lambda,\lambda _{k}})
\label{Jbound} 
\end{equation}
for  any $\bar \lambda$, $\lambda$ $\in [0,1]$, with $\bar \lambda <\lambda $, (in the  Ricci flow case, $M'$ typically is an upper bound to the scalar curvature). Let
\begin{equation}
\phi^{-1} _{\lambda _{k},\lambda}\,:\Sigma _{\lambda _{k}}\longrightarrow \Sigma _{\lambda}\;,
\end{equation}
be the inverse of $\phi _{\lambda,\lambda _{k}}$.  Denote by $d\Omega _{\lambda_{j} }$ and $d\Omega _{\lambda_{k} }$, $\lambda _{k}<\lambda _{j}$, two  probability measures  belonging to the spaces $(Prob(\Sigma _{\lambda _{j}}))$ and $(Prob(\Sigma _{\lambda _{k}}))$, respectively. Their pull--backs under the above diffeomorphism,  $(\phi^{-1} _{\lambda _{j},\lambda})^{*}\,d\Omega _{\lambda_{j} }$ and $(\phi^{-1} _{\lambda _{k},\lambda})^{*}\,d\Omega _{\lambda_{k} }$, both belong to the probability space $(Prob(\Sigma _{\lambda}))$, and we can define their Wasserstein distance at $\lambda $ according to 
\begin{eqnarray}
&& D_{2}^{W}\left((\phi^{-1} _{\lambda _{k},\lambda})^{*}\,d\Omega _{\lambda_{k} },\,(\phi^{-1} _{\lambda _{j},\lambda})^{*}\,d\Omega _{\lambda_{j} };\;\lambda \right)\doteq \label{wassertime}\\
&\doteq & \inf_{\pi_{\lambda }\in H_{\lambda } }\left( \iint_{\Sigma \times
\Sigma }d_{\lambda }(x,y)^{2}\;\pi_{\lambda } \left((\phi^{-1} _{\lambda _{k},\lambda})^{*}\,d\Omega _{\lambda_{k} },\,(\phi^{-1} _{\lambda _{j},\lambda})^{*}\,d\Omega _{\lambda_{j} }\right)\right) ^{\frac{1}{2}}\notag
\end{eqnarray}
where $H_{\lambda }\doteq H((\phi^{-1} _{\lambda _{k},\lambda})^{*}\,d\Omega _{\lambda_{k} },\,(\phi^{-1} _{\lambda _{j},\lambda})^{*}\,d\Omega _{\lambda_{j} })$ is the appropriate space of couplings, (see (\ref{wassert})), and where $d_{\lambda }(\;,\;)$ is the Riemannian distance in $(\Sigma, g(\lambda ))$. Since the cost function $d_{\lambda }(\;,\;)^{2}$ and the Jacobians $J(\phi _{\lambda,\lambda _{k}})$ are  both uniformly bounded in $\lambda $, (see (\ref{dbound}),(\ref{Jbound})), the quadratic Wasserstein distance $D_{2}^{W}(\lambda)$ defined by (\ref{wassertime}) depends smoothly from $\lambda $. In particular, we have
\begin{equation}
e^{-C''\,M''\,(\lambda -\bar \lambda)}\,D_{2}^{W}(\bar \lambda)\leq D_{2}^{W}(\lambda)\leq
e^{C''\,M''\,(\lambda -\bar \lambda)}\,D_{2}^{W}(\bar \lambda)\;,
\label{trwass}
\end{equation}
for  any $\bar \lambda$, $\lambda$ $\in [0,1]$, with $\bar \lambda <\lambda $, and for costants $C''$, $M''$ depending from the constants in (\ref{dbound}) and (\ref{Jbound}).

\noindent With these preliminary remarks along the way, let us assume that the absolutely continuous curve of probability measures 
$\varpi   :[0,1]\ni \lambda \mapsto d\varpi  _{\lambda }\in Prob(\Sigma )$, introduced in Theorem \ref{LottLength}, covers a  fiducial curve of metrics
 $\lambda \mapsto g_{ab}({\lambda })$. As above, let   $0=\lambda _{0}\leq \lambda _{1}\leq \ldots\leq \lambda _{J}=1$ be a partition of $[0,1]$. The uniform bound (\ref{dbound}) implies that we can choose $\epsilon>0$ small enough
such that
\begin{equation}
\left|d_{\lambda _{j}}(x,y)-d_{\lambda _{j-1}}(x,y) \right|\leq \left|e^{C\,M\,(\lambda _{j}-\lambda _{j-1})}-1 \right|\,d_{\lambda _{j-1}}(x,y)\;,
\label{epsi}
\end{equation}
with $\left|e^{C\,M\,(\lambda _{j}-\lambda _{j-1})}-1 \right|<<1$, whenever $\lambda _{j}-\lambda _{j-1}<\epsilon $.
Let $d\varpi _{\lambda_{j}}\in Prob(\Sigma _{\lambda _{j}})$ and $d\varpi _{\lambda_{j-1}}\in Prob(\Sigma _{\lambda _{j-1}})$ the pair of probability measures corresponding to the values $\lambda _{j}$ and $\lambda _{j-1}$ of $\lambda $. According to
(\ref{wassertime}), we can evaluate the quadratic Wasserstein distance $D_{2}^{W}\left(d\varpi _{\lambda_{j-1} },\,(\phi^{-1} _{\lambda _{j-1},\lambda})^{*}\,d\varpi _{\lambda_{j} };\;\lambda_{j-1} \right)$ between $d\varpi _{\lambda_{j}}\in Prob(\Sigma _{\lambda _{j}})$ and $d\varpi _{\lambda_{j-1}}\in Prob(\Sigma _{\lambda _{j-1}})$, at $\lambda =\lambda _{j-1}$. Thus, we define the Wasserstein length of the curve $\varpi  :[0,1]\ni \lambda \mapsto d\varpi _{\lambda }\in Prob(\Sigma )$, covering a fiducial curve of metrics $\lambda \mapsto g_{ab}(\lambda )$, according to
\begin{equation}
L_{g({\lambda })}(\varpi   )\doteq \sup_{J\in \mathbb{N}}\;\;\sup_{\lambda _{0}\leq \lambda _{1}\leq \ldots\leq \lambda _{J}}\;
\;\sum_{j=1}^{J}\; D_{2}^{W}\left(d\varpi _{\lambda_{j-1} },\,(\phi^{-1} _{\lambda _{j-1},\lambda_{j}})^{*}\,d\varpi _{\lambda_{j} };\;\lambda_{j-1} \right)   \;.
\label{walength2}
\end{equation}
Note that if $\lambda _{j-1}\leq \bar\lambda \leq \lambda _{j}$ is a refinement of the interval $[\lambda _{j-1},\,\lambda _{j}]$, then from the triangular inequality for (\ref{wassertime}) and the bound (\ref{trwass}) we get
\begin{eqnarray}
&&D_{2}^{W}\left(d\varpi _{\lambda_{j-1} },\,(\phi^{-1} _{\lambda _{j-1},\lambda_{j}})^{*}\,d\varpi _{\lambda_{j} };\;\lambda_{j-1} \right)\leq \\
&&\leq D_{2}^{W}\left(d\varpi _{\lambda_{j-1} },\,(\phi^{-1} _{\lambda _{j-1},\bar\lambda})^{*}\,d\varpi _{\bar\lambda };\;\lambda_{j-1} \right)+\notag\\
&&+D_{2}^{W}\left((\phi^{-1} _{\lambda _{j-1},\bar\lambda})^{*}d\varpi _{\bar\lambda },\,(\phi^{-1} _{\lambda _{j-1},\lambda_{j}})^{*}\,d\varpi _{\lambda_{j} };\;\lambda_{j-1} \right)\notag\\
&&\leq D_{2}^{W}\left(d\varpi _{\lambda_{j-1} },\,(\phi^{-1} _{\lambda _{j-1},\bar\lambda})^{*}\,d\varpi _{\bar\lambda };\;\lambda_{j-1} \right)+\notag\\
&&+e^{C''\,M''\,(\bar\lambda -\lambda_{j-1})}\,D_{2}^{W}\left(d\varpi _{\bar\lambda },\,(\phi^{-1} _{\bar\lambda,\lambda_{j}})^{*}\,d\varpi _{\lambda_{j} };\;\bar \lambda \right)\notag \;,
\end{eqnarray}
which implies that (\ref{walength2}) is well--behaved under refinements of the partition $0=\lambda _{0}\leq \lambda _{1}\leq \ldots\leq \lambda _{J}=1$.  

\noindent As observed above, the length of the curve $\varpi  :[0,1]\ni \lambda \mapsto d\varpi _{\lambda }\in Prob(\Sigma )$, evaluated with respect to the inner product (\ref{inner}), extends naturally to the case when the curve in question covers a fiducial curve of metrics $\lambda \mapsto g_{ab}(\lambda )$. It is sufficient to replace $|\nabla\,\Psi _{\lambda }|^{2} $in (\ref{OttoWass}) with $g^{ik}(\lambda )\,\nabla _{i}\Psi _{\lambda }\,\nabla _{k}\Psi _{\lambda }$.  
It follows that one can easily adapt the proof of proposition 3.3 in \cite{13ac} to conclude that the length of $\varpi $ with respect to the inner product (\ref{inner}) equals the length in the Wasserstein sense, \emph{i.e.},
\begin{equation}
L_{g({\lambda })}(\varpi   )=\,\int_{0}^{1}\left(\int_{\Sigma }\,g^{ik}(\lambda )\,\nabla _{i}\Psi _{\lambda }\,\nabla _{k}\Psi _{\lambda }\,d\varpi  _{\lambda } \right)^{\frac{1}{2}}d\lambda\;,
\label{ottoWasslength}
\end{equation}
which extends  the relation (\ref{OttoWass}) to the more general case considered here. 

\noindent It must be noted that a similar extension of the Hamilton--Jacobi condition (\ref{HamJac}) for characterizing Wasserstein geodesics   curves $\lambda \mapsto d\varpi_{\lambda }$, over a fiducial $\lambda \mapsto g_{ab}(\lambda )$, is quite a non--trivial problem which (to the best of my knowledge) still wait for a solution. We shall comment on this point in the concluding part of the paper and suggest a possible strategy for approaching it.

\section{Perelman's coupling for the volume-normalized Ricci flow }
To put the above probabilistic remarks in perspective we outline Perelman's characterization of the dynamics of the coupling between Ricci flow theory and scale--dependent probability measures \cite{18}.  
Let us consider the volume normalized   Ricci flow $\beta \mapsto
g_{ab}(\beta )$, $0\leq \beta <T$ \cite{6,11} associated with a metric $g_{ab}$ on a three-dimensional manifold $\Sigma $
\begin{equation}
\left\{ 
\begin{tabular}{l}
$\frac{\partial }{\partial \beta }g_{ab}(\beta )=-2{R}_{ab}(\beta )+%
\frac{2}{3}g_{ab}(\beta )\langle {R}(\beta )\rangle _{\Sigma _{\beta
}},$ \\ 
\\ 
$\;\;g_{ab}(\beta =0)=g_{ab}$,%
\end{tabular}
\right.   \label{mflow}
\end{equation}
where $R_{ab}(\beta )$ denotes the components of the Ricci tensor of $g_{ab}(\beta )$, and  
\begin{equation}
\langle {R}(\beta )\rangle _{\Sigma _{\beta }}\doteq \frac{%
\int_{\Sigma _{\beta }}{R}(\beta )  d\mu _{g(\beta )}%
}{\left[ Vol(\Sigma ,g_{ab}(\beta ))\right] }
\end{equation}
is the averaged scalar curvature with respect to the Riemannian measure $d\mu _{g(\beta )}$ defined by $g_{ab}(\beta )$.\\ One basic idea in Perelman's approach \cite{18} is to consider,
along the solution \ $g_{ab}(\beta )$ of (\ref{mflow}),  a $\beta $-dependent mapping 
\begin{gather}
f_{\beta }:\mathbb{R}\longrightarrow C^{\infty }(\Sigma_{\beta} ,\mathbb{R}) \\
\beta \longmapsto f_{\beta }:\Sigma_{\beta} \rightarrow \mathbb{R},  \notag
\end{gather}
where $C^{\infty }(\Sigma_{\beta} ,\mathbb{R})$ denotes the space of smooth
functions on $\Sigma_{\beta} $. In terms of $f_{\beta }$ one constructs on $\Sigma_{\beta} 
$\ the $\beta $-dependent measure 
\begin{equation}
d\varpi (\beta )\doteq \left( 4\pi \tau (\beta )\right) ^{-\frac{3}{2}%
}e^{-f(\beta )}d\mu _{g(\beta )},
\label{varpimeasure}
\end{equation}
where $\beta \longmapsto \tau (\beta )\in \mathbb{R}^{+}$
is a scale parameter chosen in such a way as to normalize $d\varpi (\beta )$ $\ $according to the so--called \emph{Perelman's coupling}\,:
\begin{equation}
\int_{\Sigma_{\beta} }d\varpi (\beta )=\left( 4\pi \tau (\beta )\right) ^{-\frac{3}{2%
}}\int_{\Sigma_{\beta} }e^{-f(\beta )}d\mu _{g(\beta )}=1.  \label{norm}
\end{equation}
It is easily verified that (\ref{norm}) is preserved in form along the Ricci flow (\ref{mflow}), and
\begin{equation}
\frac{d }{d \beta }\left[ \left( 4\pi \tau (\beta )\right) ^{-%
\frac{3}{2}}\int_{\Sigma_{\beta} }e^{-f(\beta )}d\mu _{g(\beta )}\right] =0,
\end{equation}
if the mapping $f_{\beta }$ and the scale parameter $\tau (\beta )$ are evolved backward in time $\beta \in (\beta ^{*},0)$ according to the coupled flows defined by 
\begin{equation}
\left\{ 
\begin{tabular}{l}
$\frac{\partial }{\partial \beta }f_{\beta }=-\Delta _{g(\beta )}f_{\beta
}+\nabla ^{i}f_{\beta }\nabla _{i}f_{\beta }-R(\beta )+\frac{3}{2}\tau
(\beta )^{-1},$  $f(\beta ^{*})=f_{0}$ \\ 
\\ 
$\frac{d }{d \beta }\tau (\beta )=\frac{2}{3}\left\langle
R(\beta )\right\rangle _{\Sigma_{\beta} (\beta )}\tau (\beta )-1,$ $\tau(\beta^{*})=\tau_{0},$%
\end{tabular}
\right.   \label{ourp}
\end{equation}
where $\Delta _{g(\beta )}$ is the Laplacian with respect to the
metric $g_{ab}(\beta )$, and $f_{0}$, $\tau_{0}$ are given (final) data, (backward $\beta $-evolution is required in order 
to have a well-posed parabolic initial value problem for (\ref{ourp})). \\

\noindent If we
assume that there exist constants (depending on $\beta $)\,\, $C^{(1)}_{\beta }$ and  $C^{(2)}_{\beta }$, $|C^{(\cdot )}_{\beta }|<\infty $, such that
\begin{equation}
C^{(1)}_{\beta }\leq \int_{0}^{\beta }\left\langle R(s)\right\rangle_{\Sigma _{s}} ds \leq  C^{(2)}_{\beta },
\label{cgrande}
\end{equation}
then
backward integration of  the $\tau (\beta )$ equation in (\ref{ourp}) along a given Ricci flow metric $\beta \longmapsto g(\beta )$, $\beta \in
\lbrack 0,T)$, provides 
\begin{equation} 
\tau (t )=e^{-\frac{2}{3}\int_{0}^{t }\left\langle R(s)\right\rangle
ds}\left[ \tau_{0}+\int_{0}^{t }e^{\frac{2}{3}\int_{0}^{\zeta
}\left\langle R(s)\right\rangle ds}d\zeta \right],
\label{soltau}
\end{equation}
where for any chosen final scale $\beta ^{\ast }\in \lbrack 0,T)$,\ we have set $t \doteq \beta ^{\ast }-\beta $ and $\tau _{0}$ $\doteq$ $\tau (t=0)$. In terms of the adimensional variable $\frac{t }{\tau _{0}}$ we can equivalently write 
\begin{equation} 
\tau (t )=\tau_{0}\,\,e^{-\frac{2}{3}\tau _{0}\int_{0}^{\frac{t}{\tau_{0}} }\left\langle R(s)\right\rangle
ds}\left[ 1+\int_{0}^{\frac{t}{\tau_{0}} }e^{\frac{2}{3}\,\tau_{0}\int_{0}^{\zeta
}\left\langle R(s)\right\rangle ds}d\zeta \right]\;.
\label{soltau2}
\end{equation}
Whereas, in terms of the forward $\frac{\beta}{\tau^{0}}$, 
\begin{equation} 
\tau (\beta )=\tau^{0}\,\,e^{\frac{2}{3}\tau ^{0}\int_{0}^{\frac{\beta }{\tau^{0}} }\left\langle R(s)\right\rangle
ds}\left[ 1-\int_{0}^{\frac{\beta}{\tau^{0}} }e^{-\frac{2}{3}\,\tau^{0}\int_{0}^{\zeta
}\left\langle R(s)\right\rangle ds}d\zeta \right]\;,
\label{soltau3}
\end{equation}
where $\tau^{0}\doteq \tau(\beta =0)$. Note that the scale parameter $\tau(t)$ is non-decreasing with $t$, along the backward Ricci flow, \emph{i.e.}\, $\frac{\partial }{\partial  t}\tau (t  )\geq 0$, as long as we  have 
$1-\frac{2}{3}\left\langle R(t  )\right\rangle _{\Sigma_{t }}\tau (t  )>0$. The geometric flow 

\begin{equation}
\left\{ 
\begin{tabular}{l}
$\frac{\partial }{\partial \beta }g_{ab}(\beta )=-2{R}_{ab}(\beta )+%
\frac{2}{3}g_{ab}(\beta )\langle {R}(\beta )\rangle _{\Sigma _{\beta
}},$ $\;\;g_{ab}(\beta =0)=g_{ab}$\\
\\
$\frac{\partial }{\partial t }f_{t }=\Delta _{g(t )}f_{t
}-\nabla ^{i}f_{t }\nabla _{i}f_{t }+R(t )-\frac{3}{2}\tau
(t )^{-1},$  $f(t=0)=f_{0}$ \\ 
\\ 
$\frac{d }{d t }\tau (t )=1-\frac{2}{3}\left\langle
R(t )\right\rangle _{\Sigma_{t} (t )}\tau (t ),$ $\tau(t=0)=\tau_{0},$%
\end{tabular}
\right.   \label{ourp2}
\end{equation}

\noindent defined by the forward volume--preserving Ricci flow (\ref{mflow}), $\beta\mapsto g_{ab}(\beta )$,\, $0\leq \beta \leq \beta ^{*}$, together  with  the backward heat and scale equations (\ref{ourp}), $t\mapsto (f_{t },\tau (t ))$,\,$t\doteq \beta ^{*}-\beta $, characterizes the \emph{Hamilton--Perelman} (volume--normalized) flow describing the coupling between the Ricci flow and the scale factorization (\ref{norm}) of the probability measure $d\varpi (t )$.
Note that one can equivalently consider the system obtained from (\ref{ourp2}) by the pull-back action, \emph{\'a la} DeTurck \cite{6b}, associated with the family of  $\beta$--dependent diffeomorphisms $\varphi :M\rightarrow M$ generated by the gradient vector field $\nabla\, f_{\beta }$. In terms of the pull--backs   $g^{*}_{ab}\doteq (\varphi ^{*}g)_{ab}$ and $f^{*}\doteq (\varphi ^{*}f)$ of the metric $g_{ab}(\beta )$ and of the function $f_{\beta }$, one can write, (see \emph{e.g.}, \cite{13e} for the detailed computation in the case of the standard Ricci flow),

\begin{equation}
\left\{ 
\begin{tabular}{l}
$\frac{\partial }{\partial \beta }g^{*}_{ab}(\beta )=-2{R}^{*}_{ab}(\beta )+2\nabla_{a} \nabla_{b} f^{*}+%
\frac{2}{3}g^{*}_{ab}(\beta )\langle {R^{*}}(\beta )\rangle _{\Sigma _{\beta
}},$ \\
\\
$\frac{\partial }{\partial t }f^{*}_{t }=\Delta _{g^{*}(t )}f^{*}_{t
}+R^{*}(t )-\frac{3}{2}(\tau^{*}
(t ))^{-1},$   \\ 
\\ 
$\frac{d }{d t }\tau^{*} (t )=1-\frac{2}{3}\left\langle
R^{*}(t )\right\rangle _{\Sigma_{t} (t )}\tau^{*} (t ),$%
\end{tabular}
\right.   \label{ourp3}
\end{equation}

\noindent where all diferential operators refer to the pull-back metric $g^{*}_{ab}$.

\noindent It must be also stressed that the equations defined by (\ref{mflow})\
and (\ref{ourp}) are based on  the standard volume preserving Ricci flow and
accordingly differ from the flows $\eta \mapsto \widetilde{g}_{ab}(\eta )$, $\eta \mapsto \widetilde{\tau }(\eta )$, and 
$\eta \mapsto f_{\eta }$, discussed by Perelman \cite{18},
\begin{equation}
\left\{ 
\begin{tabular}{l}
$\frac{\partial }{\partial \eta }\widetilde{g}_{ab}(\eta )=-2\widetilde{R}%
_{ab}(\eta ),$ \\ 
\\ 
$\widetilde{g}_{ab}(\eta =0)=\widetilde{g}_{ab},$%
\end{tabular}
\right.   \label{Perflow1}
\end{equation}
together with 
\begin{equation}
\left\{ 
\begin{tabular}{l}
$\frac{\partial }{\partial \eta }f_{\eta }=-\Delta _{\widetilde{g}(\eta
)}f_{\eta }+\nabla ^{i}f_{\eta }\nabla _{i}f_{\eta }-\widetilde{R}(\eta )+%
\frac{3}{2}\widetilde{\tau }(\eta )^{-1},$ \\ 
\\ 
$\frac{\partial }{\partial \eta }\widetilde{\tau }(\eta )=-1.$%
\end{tabular}
\right. \label{Perflow2}
\end{equation}
The flows (\ref{Perflow1}), (\ref{Perflow2}) and the ones defined by (\ref{ourp}) are
related by the usual $\eta  $-dependent homothetic rescaling \cite{11,1,6}, which maps $%
\frac{\partial }{\partial \eta }\widetilde{g}_{ab}(\eta )=-2\widetilde{R}%
_{ab}(\eta )$ to the volume-normalized Ricci flow (\ref{mflow}), \emph{i.e.}%
, 
\begin{gather}
g_{ab}(\beta (\eta ))=\left[ \frac{\int_{\Sigma }d\mu _{\widetilde{g}(\eta
=0)}}{\int_{\Sigma }d\mu _{\widetilde{g}(\eta )}}\right] ^{\frac{2}{3}}%
\widetilde{g}_{ab}(\eta ), \\
\beta (\eta )=\dint_{0}^{\eta }\left[ \frac{\int_{\Sigma }d\mu _{\widetilde{g%
}(\eta =0)}}{\int_{\Sigma }d\mu _{\widetilde{g}(s )}}\right] ^{\frac{2}{3}%
}ds,
\end{gather}
and 
\begin{equation}
\tau (\beta (\eta ))=\left[ \frac{\int_{\Sigma }d\mu _{\widetilde{g}(\eta
=0)}}{\int_{\Sigma }d\mu _{\widetilde{g}(\eta )}}\right] ^{\frac{2}{3}}%
\widetilde{\tau }(\eta ).
\end{equation}

\noindent The volume normalization makes particularly clear that, given an initial value $\tau_{0}\doteq \tau (t=0)$, the dynamics (\ref{soltau2}) of the scale $\tau (t )$ depends only on the underlying Ricci flow metric $\beta \mapsto g_{ab}(\beta )$, and not on the backward evolution of the function $f_{t}$. Thus, the localization properties of the  probability measure $d\varpi(t)$, (see (\ref{varpimeasure})) are controlled by the entropy--like quantity that one can form with $\tau (t)$ and those scalar functionals of $f_{t}$ which have the dimension of an inverse square length. The only two such objects, of geometric origin on $t\mapsto (\Sigma ,g(t))$, are the $d\varpi(t)$--expectation values of $|\nabla f_{t}|^{2}$ and of the scalar curvature $R(t)$. This latter observation immediately bring us to discuss Perelman's shrinker entropy. 

\subsection{The shrinker entropy}
The remarkable fact is that with (\ref{Perflow1}) and (\ref{Perflow2}) one can associate
Perelman's shrinker entropy $\widetilde{W}[\widetilde{g}(\eta );f_{\eta },\widetilde{\tau }(\eta )]$, defined by \cite{18} 
\begin{equation}
\widetilde{W}[\widetilde{g};f_{\eta },\widetilde{\tau }]\doteq \int_{\Sigma }\left[ \widetilde{\tau} \left( \left| \nabla
f_{\eta }\right| ^{2}+\widetilde{R}(\eta) \right) +f_{\eta }-3\right]\frac{e^{-f_{\eta }}}{(4\pi \widetilde{\tau}
(\eta ))^{\frac{3}{2}}}\,\,d\mu _{\widetilde{g}(\eta )}.
\label{W-entropyeta}
\end{equation}
The basic property of $\ \widetilde{W}[\widetilde{g}(\eta );f_{\eta },\widetilde{\tau }]$ and of the related functional 
$F[\widetilde{g};f]\doteq \int_{\Sigma }(\widetilde{R}(\eta )$\, $+|\nabla f|^{2})e^{-f}d\mu _{\widetilde{g}}$,
is that it is monotonically nondecreasing along (\ref{Perflow2}) and it provides a gradient-like structure to the Hamilton-Perelman Ricci
flow (\ref{ourp2}, \ref{ourp3}), (however, this latter interpretation must be taken with care in the probabilistic framework discussed here since, as remarked, the geometric evolution of the metric and of $f(\eta )$ are backward-conjugated).
Explicitly, one computes (see \cite{13} for a very informative analysis)
\begin{equation}
\frac{d }{d \eta }\widetilde{W}=\int_{\Sigma }2\widetilde{\tau 
}\left| \widetilde{Ric}+Hess\,f_{\eta }-\frac{1}{2\widetilde{\tau }}%
\widetilde{g}\right| ^{2}\frac{e^{-f_{\eta
}}}{(4\pi \widetilde{\tau })^{\frac{3}{2}}}\,\,d\mu _{\widetilde{g}},
\end{equation}
where $\left| ...\right| ^{2}$ is the squared $\widetilde{g}(\eta )$-norm.
Defining $\ \lambda (\widetilde{g},\widetilde{\tau })$  $\doteq $ $\
\inf_{f_{\eta }}W[\widetilde{g}(\eta );f_{\eta }]$, where the $\inf $ is
taken over all normalized $f_{\eta }$, one shows that $\inf_{\widetilde{\tau 
}>0}\lambda (\widetilde{g},\widetilde{\tau })$ is actually attained and is nondecreasing along the
Ricci flow. In particular, $\lambda (\widetilde{g},\widetilde{\tau })<0$ for
small $\widetilde{\tau }$, and $\rightarrow 0$ as $\widetilde{\tau }\searrow
0$, for any $\widetilde{g}$ on $\Sigma $. This basic property allows to probe quite
effectively the geometry of $(\Sigma ,\widetilde{g}(\eta ))$ by showing that the only shrinking Ricci flow solitons are the gradient
solitons \cite{18}.

\noindent Since the equations for $f_{\beta }$ or $f_{\eta }$ have the same  scale invariant structure, the shrinker entropy and its evolution extend, in an obvious way, to the volume normalized Hamilton-Perelman flow (\ref{ourp}), \emph{viz.}

\begin{equation}
W[g;f_{\beta },\tau ] \doteq  
\int_{\Sigma }\left[ \tau \left( \left| \nabla
f_{\beta }\right| ^{2}+R(\beta) \right) +f_{\beta }-3\right]\frac{e^{-f_{\beta }}}{(4\pi \tau
(\beta ))^{\frac{3}{2}}}\,d\mu _{g(\beta )},
\label{W-entropy}
\end{equation}
\begin{equation}
\frac{d }{d \beta }W=\int_{\Sigma }2\tau \left|
Ric+Hess\,f_{\beta }-\frac{1}{2\tau }g\right| ^{2}\frac{e^{-f_{\beta }}}{(4\pi \tau )^{\frac{3}{2}%
}}\,d\mu _{g}.
\end{equation}
What is more interesting to note is that if we introduce the normalized
Riemannian measure 
\begin{equation}
d\Pi _{\beta }\doteq Vol\left[ \Sigma_{\beta} \right] ^{-1}d\mu _{g(\beta )}
\end{equation}
associated with the (volume preserving) Ricci flow, then the $W$-functional (\ref{W-entropy}) can be be equivalently written as
\begin{gather}
W[g(\beta );f_{\beta }]\doteq \tau \,\left[I[d\varpi (\beta )\parallel d\Pi
_{\beta }] +\,\left\langle R(\beta )\right\rangle _{d\varpi (\beta )} \right]- \\
\notag \\
-S[d\varpi (\beta )\parallel d\Pi _{\beta }]+\ln \left[
Vol(\Sigma )(4\pi \tau (\beta ))^{-\frac{3}{2}}\right] -3,  \notag
\end{gather}
where
\begin{equation}
\left\langle R(\beta )\right\rangle _{d\varpi (\beta )}\doteq \int_{\Sigma
}R(\beta )d\varpi (\beta ),
\end{equation}
is the average scalar curvature with respect to $d\varpi (\beta )$, and 
\begin{equation}
I[d\varpi (\beta )\parallel d\Pi _{\beta }]\doteq \int_{\Sigma }\left|
\nabla \ln \frac{d\varpi (\beta )}{d\Pi _{\beta }}\right| ^{2}d\varpi (\beta
),
\end{equation}
\begin{equation}
S[d\varpi (\beta )\parallel d\Pi _{\beta }]\doteq \int_{\Sigma }\ln \frac{%
d\varpi (\beta )}{d\Pi _{\beta }}d\varpi (\beta ),
\end{equation}
respectively denote the entropy production functional and the relative entropy   
associated with the pair of probability measures $(d\varpi (\beta ),d\Pi
_{\beta })$. Note that the factor $Vol(\Sigma )$\,$(4\pi \tau (\beta ))^{-\frac{3}{2}}$  is generated by the normalization of $f_{\beta }\,e^{-f_{\beta }}$ 
to the probability measure density $d\varpi (\beta )/d\Pi _{\beta }$, and one may equivalently write
\begin{eqnarray}
&&S[d\varpi (\beta )\parallel d\Pi _{\beta }]-\ln \left[
Vol(\Sigma )(4\pi \tau (\beta ))^{-\frac{3}{2}}\right]=\\
&=&-Vol(\Sigma )(4\pi \tau (\beta ))^{-\frac{3}{2}}\,\int_{\Sigma }f_{\beta }\,e^{-f_{\beta }}\,d\Pi _{\beta }\;\notag.
\end{eqnarray}
Note also that we can rewrite  $W[g(\beta
);f_{\beta }]$ in terms of  the defective LSI  functional (\ref{logsob}), \emph{i.e.}
\begin{equation}
W[g(\beta );f_{\beta }]= LSI[(2\tau (\beta ))^{-1};\,\left\langle R(\beta )\right\rangle _{d\varpi (\beta )}]+
\ln \left[\frac{
Vol(\Sigma )}{(4\pi \tau (\beta ))^{\frac{3}{2}}}\right] -3\;.
\label{defectiveW}  
\end{equation}
Whereas this is  rather trivial consequence of the fact that  the standard Gaussian logarithmic Sobolev inequality \cite{10,18} lies at the origin of the definition (\ref{W-entropyeta}) of $W[g(\beta );f_{\beta }]$, it also indicates that $W[g(\beta );f_{\beta }]$ is not an entropy, but rather an entropy--balance functional controlling the rate of variation of $S[d\varpi (\beta )\parallel d\Pi _{\beta }]$ along the Ricci flow, and hence the localization properties of $d\varpi (\beta )$. To prove this latter remark, let us consider, as in the above analysis, $\beta \longmapsto g_{ab}(\beta )$, $\beta \in
\lbrack 0,T)$, \ with $t\doteq \beta ^{\ast }-\beta $, $\beta ^{\ast }\in \lbrack 0,T)$. In the backward direction, along
the given volume-normalized Ricci flow, we have 
\begin{equation}
\frac{\partial }{\partial t}d\Pi _{t}=\left[R(t)- \left\langle R(t
)\right\rangle_{\Sigma _{t}} \right] d\Pi _{t },
\label{bdet}
\end{equation}
moreover, Perelman's condition (\ref{ourp})  yields the conjugate heat equation
\begin{equation}
\frac{\partial }{\partial t }\left[ \frac{d\varpi (t )}{d\Pi _{t
}}\right] =\Delta _{g(t )}\left[ \frac{d\varpi (t )}{d\Pi _{t }}%
\right] -\frac{d\varpi (t )}{d\Pi _{t }}\left[ R(t
)-\left\langle R(t )\right\rangle_{\Sigma _{t}} \right] ,  \label{lindensity}
\end{equation}
where 
\begin{equation}
\frac{d\varpi (t )}{d\Pi _{t }}=\left( 4\pi \tau (t )\right) ^{-%
\frac{3}{2}}Vol[\Sigma_{t} ]e^{-f(t )}.
\end{equation}
Since $d\Pi _{t }$ is covariantly constant with respect to the Levi-Civita connection $\nabla $ associated with $g(t)$, 
we can exploit (\ref{bdet}) and write (\ref{lindensity}) as the (non-uniformly parabolic) probability diffusion  $(d\varpi (t))_{t\geq 0}$ PDE,
\begin{equation}
\left\{ 
\begin{tabular}{l}
$\frac{\partial }{\partial t}d\varpi (t)=\Delta _{g(t)}\left( d\varpi
(t)\right) , \,\, t\doteq \beta ^{\ast }-\beta 
$ \\ 
\\ 
$d\varpi (t=0)=d\varpi_{0}.$%
\end{tabular}
\right.   \label{heat}
\end{equation}
From (\ref{bdet}) and (\ref{heat}),  we get
\begin{eqnarray}
&&\frac{d}{d\,t }S[d\varpi (t )\parallel d\Pi _{t }]=\label{SIRR}\\
\notag\\
&=&\int_{\Sigma }\ln\,\frac{d\varpi (t )}{d\Pi _{t }}\,\Delta _{g(t)}\left( d\varpi
(t)\right)-\int_{\Sigma }\left[R(t)- \left\langle R(t
)\right\rangle \right]\,d\varpi (t )\notag\\
\notag\\
&=&-I[d\varpi (t )\parallel d\Pi
_{t}] -\,\left\langle R(t )\right\rangle _{d\varpi (t )} +
\,\,\left\langle R(t )\right\rangle _{\Sigma_{t} }\notag\;.
\end{eqnarray}

\noindent On the other hand, according to the evolution (\ref{ourp2}) of the scale parameter $\tau (t)$, we can write
\begin{eqnarray}
\left\langle R(t )\right\rangle _{\Sigma_{t} }&=&\frac{3}{2}\,\tau^{-1} (t)-\frac{3}{2}\,\tau^{-1} (t)\,\frac{d }{d t}\,\tau (t)\\
&=& \frac{3}{2}\,\tau^{-1} (t)-\,\frac{d }{d t}\,\ln \frac{(4\pi\, \tau (t))^{\frac{3}{2}}}{Vol(\Sigma )}\notag\;,
\label{taudim}
\end{eqnarray}
where, in the last line, we have normalized $(4\pi\, \tau (t))^{\frac{3}{2}}$ to $Vol(\Sigma )$ for dimensional reasons, (and we have introduced the factor $4\pi $ for later convenience).  Inserting (\ref{taudim}) in (\ref{SIRR}), we get
\begin{eqnarray}
&&\frac{d }{d t}\,\left(S[d\varpi (t )\parallel d\Pi _{t }]+\ln 
\frac{(4\pi \tau (t ))^{\frac{3}{2}}}{Vol(\Sigma )}\right)\label{Wflow} \\
\notag\\
&=&-\,\left[ I[d\varpi (t )\parallel d\Pi
_{t}] +\,\left\langle R(t )\right\rangle _{d\varpi (t )}  \right] + \frac{3}{2}\,\tau (t )^{-1}\notag\\
\notag\\
&=&-(4\pi \tau (t ))^{-\frac{3}{2}}\,F[g,f]\,+\,\frac{3}{2}\,\tau (t )^{-1}\notag\;,
\end{eqnarray}
where 
\begin{equation}
F[g,f]\doteq \int_{\Sigma }({R}(t )\,+|\nabla f|^{2})e^{-f}d\mu _{{g}}\;,
\end{equation}
is the standard Perelman functional associated with the volume--preserving flow  (\ref{ourp2}). (I wish to thank  Alessio Figalli for suggesting that a relation of this type should hold; the fact that the entropy functional $F[g;f]$, for the standard un--normalized, Ricci flow, is the time--derivative of the entropy--like quantity $\int_{\Sigma }f_{t }\,e^{-f_{t }}\,d\mu _{g(t) }$ has been noticed also in \cite{13e}). The relation (\ref{Wflow})
characterizes also the shrinker entropy $W[g(\beta );f_{\beta }]$ as the variation, along the backward Ricci flow, of the 
functional
\begin{eqnarray}
&&\mathcal{G}[d\varpi (t ),\, d\Pi _{t },\,\tau (t )]\doteq   S[d\varpi (t )\parallel d\Pi _{t }]+\ln 
\frac{(4\pi \tau (t ))^{\frac{3}{2}}}{Vol(\Sigma )}+\frac{3}{2}\\
&=&-Vol(\Sigma )(4\pi \tau (t ))^{-\frac{3}{2}}\,\int_{\Sigma }f_{t }\,e^{-f_{t }}\,d\Pi _{t }+\frac{3}{2}\notag\;.
\end{eqnarray}
Indeed, (\ref{Wflow}) can be rewritten as 
\begin{equation}
\frac{d}{d\,t}\left( \tau (t)\mathcal{G}\right)
 =-W[g(\beta );f_{\beta }]-\frac{2}{3}\,\left\langle R(t )\right\rangle _{\Sigma_{t} }\,\tau (t)
\mathcal{G}\;.
\label{wasder}
\end{equation}

\noindent  The fact that 
(\ref{wasder}) has the same formal structure of the evolution equation (\ref{ourp2}) of the scale parameter $\tau (t)$ further confirms that $W[g(\beta );f_{\beta }]$ controls, via the effective scale 
$\tau_{eff} (t)\doteq \tau (t)\mathcal{G} $, the localization properties of the measure $d\varpi (t)$.

\noindent It is important  to remark that, according to the relation (\ref{SIRR}), the relative entropy 
$S[d\varpi (t )\parallel d\Pi _{t}]$ is not monotonic along the Hamilton--Perelman flow. This is related to the fact that the evolution $t\mapsto \frac{d\varpi (t)}{d\Pi _{t}}$, described by (\ref{lindensity}), is not a gradient flow for 
$S[d\varpi (t )\parallel d\Pi _{t}]$. The deviation from being gradient--like are due to the presence of the curvature fluctuation term $R(t)-\left\langle R(t )\right\rangle_{\Sigma _{t}}$ in (\ref{lindensity}), (under $d\varpi (t)$--expectation, this term yields $\left\langle R(t )\right\rangle_{d\varpi {t}}-\left\langle R(t )\right\rangle_{\Sigma _{t}}$). According to the analysis presented in section \ref{Ottop}, a possible strategy for compensating such a fluctuation term is to trade $R(t)-\left\langle R(t )\right\rangle_{\Sigma _{t}}$  for a ($t$--dependent) potential, so as to transform the heat diffusion $(d\varpi (t))_{t\geq 0}$ into a Fokker--Planck diffusion.  However, before discussing such a transformation in detail, there is still an elementary but interesting property of (\ref{SIRR}) that we would like to point out and concerning the role and the monotonicity properties of the $d\varpi (t)$ average $\left\langle R(t )\right\rangle_{d\varpi {t}}$ of the scalar curvature.

\subsection{A renormalized curvature entropy}
 To begin with, let us observe that the basic relation (\ref{Wflow}) becomes particularly simple if we rescale the parameter $\tau (t)$ according to
\begin{equation}
\hat {\tau }(t)\doteq \tau_{0}\,\,e^{-\frac{2}{3}\tau_{0}\int_{0}^{\frac{t}{\tau_{0}}}\,\left\langle R(s)\right\rangle _{\Sigma
_{s}}ds}=
\frac{\tau (t)}{1+\int_{0}^{\frac{t}{\tau_{0}}}e^{\frac{2}{3}\tau_{0}\int_{0}^{\zeta
}\left\langle R(s)\right\rangle ds}d\zeta}\;.
\label{rentau}
\end{equation}
Note that $\hat {\tau }(t)$ satisfies
\begin{equation}
\frac{d}{dt}\hat {\tau }(t)=-\frac{2}{3}\,\left\langle R(t)\right\rangle _{\Sigma
_{s}}\,\hat {\tau }(t)\;.
\end{equation}

\noindent It is easily verified that in terms of $\hat {\tau }(t)$ we can rewrite (\ref{Wflow}) as 

\begin{eqnarray}
&&\frac{d }{d t}\,\left(S[d\varpi (t )\parallel d\Pi _{t }]+\ln 
\frac{(4\pi \hat{\tau} (t ))^{\frac{3}{2}}}{Vol(\Sigma )}\right) \label{RenWflow}\\
\notag\\
&=&-\,\left[ I[d\varpi (t )\parallel d\Pi
_{t}] +\,\left\langle R(t )\right\rangle _{d\varpi (t )}  \right] \;.\notag
\end{eqnarray}
Since 
\begin{eqnarray}
&&S[d\varpi (t )\parallel d\Pi _{t }]+\ln 
\frac{(4\pi \hat{\tau} (t ))^{\frac{3}{2}}}{Vol(\Sigma )}\\
\notag\\
&=&\int_{\Sigma }d\varpi (t )\,\ln\left[ \left(
\frac{(4\pi \hat{\tau} (t ))^{\frac{3}{2}}}{Vol(\Sigma )}\right)\frac{d\varpi (t )}{d\Pi _{t }}  \right]\,,\notag
\end{eqnarray}
and 
\begin{equation}
\int_{\Sigma }d\varpi (t )\,\left|\nabla \ln\left[ \left(
\frac{(4\pi \hat{\tau} (t ))^{\frac{3}{2}}}{Vol(\Sigma )}\right)\frac{d\varpi (t )}{d\Pi _{t }}  \right]\right|^{2}=I[d\varpi (t )\parallel d\Pi_{t}]\;,
\end{equation}

\noindent it follows that $\left\langle R(t )\right\rangle _{d\varpi (t )}$ represents the obstruction to $(d\varpi (t))_{t\geq 0}$ for being a gradient--like flow. It is indeed  the LSI defective parameter 
$\left\langle R\right\rangle _{d\varpi }$ appearing in the shrinker entropy.
One can easily check that $\left\langle R(t )\right\rangle _{d\varpi (t )}$ is not monotonic along the backward Ricci flow,  however $\hat{\tau }(t)\left\langle R(t )\right\rangle _{d\varpi (t )}$ turns out to be weakly--monotonic, and we have the following
\begin{theorem}
For a given Ricci flow metric $\beta \longmapsto g(\beta )$, $\beta \in
\lbrack 0,T)$, and for any chosen $\beta ^{\ast }\in \lbrack 0,T)$,\ let $%
t\longmapsto d\varpi (t)$, and $t\longmapsto \tau (t)$, $t\doteq \beta ^{\ast }-\beta $, be the solutions of
\ (\ref{heat}) and (\ref{soltau2}) corresponding to the initial data $%
d\varpi _{0}$ and $\tau_{0}$, respectively. Then, along the backward Ricci flow, we have

\begin{equation}
\frac{\partial }{\partial t}\left[ \hat{\tau }(t)\left\langle R(t)\right\rangle _{d\varpi (t)}
 \right] =-2\hat{\tau }(t)\int_{\Sigma }d\varpi (t)\left| {Ric}(t)\right| ^{2}\leq 0.
\end{equation}
In particular, 

\begin{equation}
\tau_{0}\,\,e^{-\frac{2}{3}\tau_{0}\int_{0}^{\frac{t}{\tau_{0}}}\,\left\langle R(s)\right\rangle _{\Sigma
_{s}}ds}
\left\langle R(t)\right\rangle _{d\varpi (t)}
\label{monot}
\end{equation}

\noindent is nonincreasing as a function of $t\in [0,\beta ^{*})$.
\label{curvtau}
\end{theorem}

\begin{proof} Let us recall that along the (volume-normalized)
Ricci flow $\beta \mapsto g_{ik}(\beta )$, $\beta \in (0,T)$ we have \cite{6}  
\begin{equation}
\frac{\partial \,R(\beta )}{\partial \beta }=\Delta _{g(\beta )}R(\beta )+2|%
\widehat{Ric}(\beta )|^{2}+\frac{2}{3}R(\beta )(R(\beta )-\left\langle
R(\beta )\right\rangle _{\Sigma _{\beta }}),
\label{evscalcurv}
\end{equation}
where $\widehat{Ric}$ denotes the trace-free part $R_{ab}-\frac{1}{3}g_{ab}R$ of the Ricci tensor.  A  direct computation exploiting (\ref{evscalcurv}) provides
\begin{gather}
\frac{\partial }{\partial t}\int_{\Sigma }d\varpi (t)\hat{\tau }(t)R(t)=-\frac{2}{3}%
\hat{\tau }(t)\langle R(t)\rangle _{\Sigma _{t}}\int_{\Sigma }d\varpi (t)R(t)+ \\
-\hat{\tau }(t)\int_{\Sigma }\Delta _{g(t)}R(t)d\varpi (t)-2\hat{\tau }(t)\int_{\Sigma }d\varpi
(t)\left| \widehat{Ric}(t)\right| ^{2}-  \notag \\
-\frac{2}{3}\hat{\tau }(t)\int_{\Sigma }d\varpi (t)R(t)\left( R(t)-\langle R(t)\rangle
_{\Sigma _{t}}\right) +  \notag \\
+\hat{\tau }(t)\int_{\Sigma }R(t)\Delta _{g(t)}(d\varpi (t)).  \notag
\end{gather}
Since $\int_{\Sigma }\Delta _{g(t)}R(t)d\varpi (t)$ $=$
$\int_{\Sigma }R(t)\Delta _{g(t)}(d\varpi (t))$ and $|\widehat{Ric}|^{2}+\frac{1}{3}R(t)^{2}$ $=$
$|Ric|^{2}$, we get the stated result.
\end{proof}

\noindent The monotonicity of (\ref{monot}) immediately implies that
\begin{equation}
\left\langle R(t_{1})\right\rangle _{d\varpi (t_{1})}\geq
\left\langle R(t_{2})\right\rangle _{d\varpi (t_{2})}\;
 \exp\left[{-\frac{2}{3}\tau_{0}%
\int_{{t_{1}}\backslash {\tau_{0}}}^{{t_{2}}\backslash {\tau_{0}}}\;\left\langle R(s)\right\rangle _{\Sigma
_{s}}ds}\right]\;,
\end{equation}
for any $t_{2}\geq t_{1}$. Thus, if $\left\langle R(t_{1})\right\rangle
_{d\varpi (t_{1})}\leq 0$, the measure $\{d\varpi (t)\}_{t\geq t_{1}}$ will
diffuse in a region of \ $(\Sigma ,g(t))$ where  $R(t)\leq 0$, (regardless
of the sign of the overall average $\left\langle R(t)\right\rangle _{\Sigma
_{t}}$). Since $t$ parametrizes the backward flow, this remark implies that
regions of negative scalar curvature result from the Ricci flow $\beta $%
-evolution of regions with negative scalar curvature (in other words
localized positive scalar curvature cannot evolve into negative scalar
curvature under the Ricci flow).  \ This is another manifestation of the
fact that the flow prefers positive scalar curvature.\

\noindent More generally, we can rewrite the evolution equation in Theorem \ref{curvtau} as
\begin{equation}
\frac{d }{d t}\left[ \hat{\tau }(t)\left\langle R(t)\right\rangle
_{d\varpi (t)}\right] =-2\hat{\tau }(t)\left\langle \left| Ric(t)\right|
^{2}\right\rangle _{d\varpi (t)}
\end{equation}
which, by factorizing the Ricci tensor in its trace-free and trace part, and by adding
and subtracting the term $\frac{2}{3}\widetilde{\tau }(t)\left\langle R(t)\right\rangle
_{d\varpi (t)}^{2}$, yields
\begin{gather}
\frac{d }{d t}\left[ \hat{\tau }(t)\left\langle R(t)\right\rangle
_{d\varpi (t)}\right] =-2\hat{\tau }(t)\left\langle \left| \widehat{Ric(t)}\right|
^{2}\right\rangle _{d\varpi (t)}- \\
-\frac{2}{3}\hat{\tau }(t)\left( \left\langle R^{2}(t)\right\rangle _{d\varpi
(t)}-\left\langle R(t)\right\rangle _{d\varpi (t)}^{2}\right) -\frac{2}{3}%
\hat{\tau }(t)\left\langle R(t)\right\rangle _{d\varpi (t)}^{2}.  \notag
\end{gather}
Since the terms on the right hand side of this expression are all
non-negative, we get (after dividing and multiplying by $\hat{\tau }(t)>0$) 
\begin{equation}
\frac{d }{d t}\left[ \hat{\tau }(t)\left\langle R(t)\right\rangle
_{d\varpi (t)}\right] \leq -\frac{2}{3}\hat{\tau }^{-1}(t)\left[ \hat{\tau }(t)\left\langle
R(t)\right\rangle _{d\varpi (t)}\right] ^{2}.
\end{equation}
This integrates to the Harnack--type inequality
\begin{equation}
\left\langle R(t)\right\rangle _{d\varpi(t)}\leq \frac{e^{\frac{2}{3}\tau _{0}%
\int_{0}^{\frac{t}{\tau _{0}}}\;\left\langle R(s)\right\rangle ds}\left\langle R_{0}\right\rangle
_{d\varpi_{0}}}{1+\frac{2}{3}\tau _{0}\left\langle R_{0}\right\rangle _{d\varpi_{0}
}\;\int_{0}^{\frac{t}{\tau _{0}}}e^{\frac{2}{3}\tau _{0}\int_{0}^{\zeta }\left\langle
R(s)\right\rangle ds}d\zeta },
\end{equation}
where $0\leq \zeta \leq t$, and $R_{0}\doteq{R(t=0)}$. Since $\left\langle R_{0}\right\rangle _{d\varpi _{0}}$ results from the $\beta $-evolution 
of $R$ in the regions localized by the measure $d\varpi (t)$, the above estimate allows, as we have seen above, to compare scalar curvature at different times in different regions along the backward Ricci flow.

\section{Ricci flow and Fokker-Planck diffusion on $Prob(\Sigma)$}
\noindent 
According to (\ref{bdet}) and (\ref{SIRR}), the curvature fluctuation term $R(t)-\left\langle R(t )\right\rangle_{\Sigma _{t}}$ drives the dynamics of the Riemannian measure $d\Pi _{t}$, and obstructs the gradient--like nature of Perelman's diffusion $(d\varpi (t))_{t\geq 0}$. Such a  behavior is in line with the geometric properties of Otto's parametrization of diffusion processes, along a fiducial curve of metrics, discussed in section \ref{Ottop}. In this connection, let us recall that the curve in the space of Riemannian metrics
\begin{eqnarray}
\mathcal{R}iem(\Sigma ) &\longrightarrow& \mathcal{R}iem(\Sigma )\\
(\Sigma ,g) &\mapsto& (\Sigma ,g(\beta ))\;, \nonumber
\end{eqnarray}
defined by the Ricci flow (\ref{mflow}),  is natural in the geometrical sense since it is $\mathcal{D}iff(\Sigma )$--equivariant, and  always admits a solution $\beta \mapsto (\Sigma ,g(\beta ))$, in a maximal interval $0\leq \beta \leq T_{0}$,  for some $T_{0}\leq \infty $, (if such a $T_{0}$ is finite, the flow 
necessarily develops \cite{11} a curvature singularity as $\beta \nearrow T_{0}$, {\emph{i.e.}}, 
$\lim_{\beta \nearrow T_{0}}\, [\sup_{x\in \Sigma }\,|Rm(x,\beta )|]=\infty $, where $Rm(\beta )\doteq (R^{l}_{ijk}(\beta ))$ is the Riemann tensor of $(\Sigma ,g(\beta ))$). Such a geometrical naturality is the basic reason why, as we shall see, the Ricci flow  provides  a natural class of fiducial curve of metrics along which Otto's parametrization turns out to be particularly effective.

\subsection{A potential for scalar curvature fluctuations}
\label{potentialsect}

Let us start by observing that the curve of Riemannian measures $t\mapsto d\Pi _{t}$ can be formally considered as the lift
to the bundle $Prob(\Sigma )$ of the curve in $\mathcal{R}iem(\Sigma )$ defined by the  backward Ricci flow $[0,\beta ^{\ast }]\ni t\mapsto g_{ab}(t)$. In this way, one characterizes a fiducial curve of reference measures in the bundle 
$Prob(\Sigma )$, and for each $t\in [0,\beta ^{*}]$ we can naturally describe the corresponding fiber of $Prob(\Sigma )$
 over $g_{ab}(t)$, as
\begin{equation}
Prob(\Sigma_{t} )\doteq Prob(\Sigma,g(t) )= \left\{N\,d\Pi_{t}: N\in C_{b}(\Sigma ,\mathbb{R}^{+}),\,\int_{%
\Sigma }N\,d\Pi_{t} =1\right\}\;.
\end{equation} 
To fully exploit such a description we need to characterize  the tangent vector  to 
the fiducial curve $t\rightarrow d\Pi _{t}$. Let us consider a generic value of the parameter $t$, say $t=s$.\ Since 
$\int_{\Sigma_{t} }\frac{%
\partial }{\partial t}d\Pi _{t}|_{t=s}=0$, we have that $\frac{%
\partial }{\partial t}d\Pi _{t}|_{t=s} \in T_{d\Pi _{t}}Prob(\Sigma )$.\ According to (\ref{otto}), and in analogy with (\ref{elliptictheta}) we  parametrize  $\frac{%
\partial }{\partial t}d\Pi _{t}|_{t=s}$ in terms of a scalar curvature-fluctuations potential $\Phi
_{s}$ obtained as the solution of the elliptic equation 
\begin{equation}
g^{ik}(s)\nabla _{i}\left( d\Pi _{s}\,\nabla _{k}\Phi _{s}\right) =-\left. 
\frac{\partial }{\partial t}d\Pi _{t}\right| _{t=s},
\label{ellipticR1}
\end{equation}
where $\frac{\partial }{\partial t}d\Pi _{t}|_{t=s}$ is given by (\ref{bdet}) for each given $t=s$, {\em i.e.}, (again exploiting the covariant
constancy of $d\Pi _{t}$),
\begin{equation}
\Delta _{g(s)}\Phi _{s}=-\left( R(s)-\left\langle R(s)\right\rangle \right) .
\label{ellipticR}
\end{equation}
Formally, given a solution $\Phi _{s}$ of (\ref{ellipticR}) and any ($t$-independent) smooth function with compact support $\zeta \in C_{0}^{\infty }(\Sigma, \mathbb{R})$, we can write 
\begin{equation}
\frac{d }{d t}\int_{\Sigma_{t} }\zeta d\Pi _{t}=\int_{\Sigma_{t} }\zeta \frac{\partial }{\partial t}d\Pi _{t}
=\int_{\Sigma_{t}
}\nabla ^{i}\zeta \,\nabla _{i}\Phi _{t}\,d\Pi _{t}= \left\langle \zeta ,\Phi _{t}\right\rangle_{d\Pi _{t}}.
\label{fi} 
\end{equation}
According to (\ref{inner}) the 
relation (\ref{fi}) identifies $d\Pi _{t}\longmapsto \Phi _{t} $
as the tangent vector to the curve $t\rightarrow d\Pi _{t}$ and defines the curvature (fluctuation) potential in which a probability density, evolving along a Ricci flow manifold, diffuses.
 
\noindent Since $\int_{\Sigma }d\Pi _{s}\left( R(s)-\left\langle R(s)\right\rangle\right) =0$, and
\begin{equation}
\int_{\Sigma _{t}}(R(t)-\langle R(t)\rangle _{\Sigma _{t}})\Phi _{t}\,\,d\Pi_{t}=\int_{\Sigma _{t}}\left| \nabla \,\Phi _{t}\right| ^{2}\,d\Pi _{t}
\end{equation}
we get that, as long as the
(volume-normalized) Ricci flow exists, equation (\ref{ellipticR}) admits a solution unique
up to constants. The $L^{2}((\Sigma ,d\Pi _{s}),\mathbb{R})$ norm
of $\left( R(s)-\left\langle R(s)\right\rangle \right) $ is given by 
$\left\langle R^{2}(s)\right\rangle -\left\langle
R(s)\right\rangle ^{2}$, thus $\Phi _{s}$ is in the Sobolev
space $H_{2}((\Sigma ,d\Pi _{s}),\mathbb{R})$ if the mean square fluctuations in the scalar curvature are bounded.
More generally, we know (see \emph{e.g.} \cite{6}) that if $%
\beta \mapsto g_{ab}(\beta )$ is a solution of the Ricci flow equation (for
which the weak maximum principle holds), then bounds on the curvature (and its derivatives) of the
initial metric induce a priori bounds on  all derivatives $|\nabla
^{m}R(x,s)|$ for a sufficiently short time. Thus, for any given $t=s$ for which the Ricci flow is non-singular, we can
assume that $\left( R(s)-\left\langle R(s)\right\rangle \right) $ is $%
C^{\infty }(\Sigma ,\mathbb{R})$ and by elliptic regularity we get that $\Phi _{s}$ $\in$ 
$C^{\infty }(\Sigma ,\mathbb{R})$. Along the same lines, we also have 
\begin{lemma}
For any  $\Phi _{t}\in C^{3}(\Sigma, \mathbb{R})$ solution of \ (\ref{ellipticR}) the following relation holds
\begin{gather}
\int_{\Sigma _{t}}\left| Hess\,\Phi _{t}\right| ^{2}\,d\Pi _{t}+\int_{\Sigma
_{t}}R^{ik}(t)\,\nabla _{i}\Phi _{t}\nabla _{k}\Phi _{t}\,d\Pi _{t}= \label{interesse} \\ 
=\langle R(t)^{2}\rangle _{\Sigma _{t}}-\langle R(t)\rangle _{\Sigma
_{t}}^{2}.  \notag
\end{gather}
If the Ricci curvature of $(\Sigma _{t},g(t))$ is positive,  then 
\begin{equation}
\int_{\Sigma _{t}}\left| \nabla \,\Phi _{t}\right| ^{2}\,d\Pi _{t}\,\leq 
\frac{2}{3}\frac{\langle R(t)^{2}\rangle _{\Sigma _{t}}-\langle R(t)\rangle
_{\Sigma _{t}}^{2}}{K_{t}},
\label{positive}
\end{equation}
where $K_{t}>0$ is the lower bound of $Ric(t)$. Moreover, if $L_{g(t)}(d\Pi _{t})$
denotes the (quadratic) Wasserstein length of the curve $t\mapsto d\Pi _{t}$, then we have
\begin{equation}
L_{g(t)}(d\Pi _{t})\leq \sqrt{\frac{2}{3}}\int_{0}^{\beta ^{*}}\left( \frac{\langle R(t)^{2}\rangle _{\Sigma _{t}}-\langle R(t)\rangle
_{\Sigma _{t}}^{2}}{K_{t}}  \right)^{\frac{1}{2}}\,dt\;.
\end{equation}
\end{lemma}
\begin{proof}
These results are elementary consequences of the Ricci commutation relation $\nabla ^{j}\nabla _{i}\nabla _{j}\,f$\,\,$-\nabla _{i}\nabla ^{j}\nabla _{j}\,f$\,$=R_{ij}\nabla ^{j}\,f$, valid for any $f\in C_{3}(\Sigma, \mathbb{R})$. In particular, for $\Phi _{t}\in C_{3}(\Sigma,\mathbb{R})$, 
consider the expression  $2\nabla ^{i}\Phi _{t}\nabla _{i}\Delta _{g(t)}\Phi
_{t}$. The Bochner-Weitzenb\"{o}ck formula provides
\begin{equation}
2\nabla ^{i}\Phi _{t}\nabla _{i}\Delta _{g(t)}\Phi _{t}=\Delta _{g(t)}\left|
\nabla \Phi _{t}\right| ^{2}-2\left| Hess\,\Phi _{t}\right|
^{2}-2R^{ik}(t)\,\nabla _{i}\Phi _{t}\nabla _{k}\Phi _{t},
\end{equation}
pointwise. Thus
\begin{gather}
-\int_{\Sigma _{t}}\left| Hess\,\Phi _{t}\right| ^{2}\,d\Pi
_{t}-\int_{\Sigma _{t}}R^{ik}(t)\,\nabla _{i}\Phi _{t}\nabla _{k}\Phi
_{t}\,d\Pi _{t}= \\
=\int_{\Sigma _{t}}\nabla ^{i}\Phi _{t}\nabla _{i}\Delta _{g(t)}\Phi
_{t}\,\,d\Pi _{t}=\int_{\Sigma _{t}}(R(t)-\langle R(t)\rangle _{\Sigma
_{t}})\Delta _{g(t)}\Phi _{t}\,\,d\Pi _{t}=  \notag \\
=-\left( \langle R(t)^{2}\rangle _{\Sigma _{t}}-\langle R(t)\rangle _{\Sigma
_{t}}^{2} \right),  \notag
\end{gather}
where we have integrated by parts and exploited (\ref{ellipticR}). 
Since 
\begin{equation}
\int_{\Sigma _{t}}\left| Hess\,\Phi _{t}\right| ^{2}\,d\Pi _{t}\geq \frac{1}{%
3}\int_{\Sigma _{t}}\left| \Delta _{g(t)}\,\Phi _{t}\right| ^{2}\,d\Pi _{t}=%
\frac{1}{3}(\langle R^{2}\rangle _{\Sigma _{t}}-\langle R\rangle _{\Sigma
_{t}}^{2}),
\end{equation}
we get from (\ref{interesse})
\begin{equation}
\frac{2}{3}\left( \langle R(t)^{2}\rangle _{\Sigma _{t}}-\langle R(t)\rangle
_{\Sigma _{t}}^{2}\right) \geq \int_{\Sigma _{t}}R^{ik}(t)\,\nabla _{i}\Phi
_{t}\nabla _{k}\Phi _{t}\,d\Pi _{t},
\end{equation}
which, if the Ricci curvature has a positive lower bound $K_{t}$, yields (\ref{positive}).
Note that, since positive Ricci curvature is preserved along the Ricci
flow, the bound (\ref{positive}) holds for every $t$. If $L_{g(t)}(d\Pi _{t})$
is the (quadratic) Wasserstein length of the curve $t\mapsto d\Pi _{t}$, defined according to (\ref{walength2}), then the identification (\ref{ottoWasslength}) and the bound (\ref{positive})  imply 
\begin{eqnarray}
L_{g(t)}(d\Pi _{t})&=&\int_{0}^{\beta ^{*}}\left(\int_{\Sigma _{t}}\left| \nabla \,\Phi _{t}\right| ^{2}\,d\Pi _{t}  \right)^{\frac{1}{2}}\,dt\\
&\leq& \sqrt{\frac{2}{3}}\int_{0}^{\beta ^{*}}\left( \frac{\langle R(t)^{2}\rangle _{\Sigma _{t}}-\langle R(t)\rangle
_{\Sigma _{t}}^{2}}{K_{t}}  \right)^{\frac{1}{2}}\,dt\;,\notag
\end{eqnarray}
as stated.
\end{proof}

\noindent Note that equation (\ref{ellipticR}) has a familiar counterpart in the Ricci flow theory for surfaces \cite{11a}, (see also \cite{6}). It arises also in 
K\" ahler geometry where it provides the relation between the K\" ahler Ricci potential and the scalar curvature.  Moreover, if we consider a gradient Ricci soliton $Ric(g(t))=Hess \, \Psi_{t}$, 
\emph{i.e.}, a fixed point of the flow obtained by quotienting
the space of metrics under diffeomorphisms
and scalings, then clearly $\Phi _{t}\doteq -\,\Psi _{t}$. In the  general case, (\emph{viz.}, when $(\Sigma ,g)$ is neither K\" ahler or a gradient soliton), we can still obtain a geometrical characterization of $\Phi _{s}$. We start by deriving an asymptotic expression for $\Phi _{s}$ valid in local geodesic
coordinates (LGC), $\{x^{i}\}_{LGC}$, at any given fixed point $p\in \Sigma _{s}$. Let us denote by $r(x)\doteq d(p,x)$ 
the Lipschitz function providing the distance from $p$ to $x$. For $x\notin Cut(p)$, the cut locus of $p$, we set 
$x^{j}=r\,u^{j}$, with $u^{j}$ coordinates on the unit sphere $\mathbb{S}^{2}\subset T_{p}\Sigma_{s} $. The pull-back of the Riemannian measure $d\Pi_{s} $ under the exponential mapping $exp_{p}:\,T_{p}\Sigma_{s} \rightarrow \Sigma_{s} $, provides the familiar asymptotics in geodesic polar coordinates
\begin{equation}
\frac{exp_{p}^{*}(d\Pi _{s})}{d\Pi _{Eucl.}}
\underset{LGC}{=}1-\frac{1}{6}\,R_{ik}(p)\,x^{i}x^{k}-\frac{1}{12}%
\,\nabla _{j}\,R_{ik}(p)\,x^{i}x^{k}x^{j}+O(r^{4}),
\end{equation}
where $d\Pi _{Eucl.}$ is the standard Euclidean volume element in polar coordinates in $T_{p}\Sigma _{s}$, and $R_{ik}(p)$ are the components of the Ricci tensor at $p
$. Thus, if we take the (Euclidean) Hessian of the function \,\, $6\left(\frac{exp_{p}^{*}(d\Pi _{s})}{d\Pi _{Eucl.}}\right)$ \, 
$+\frac{1%
}{3}\left\langle R\right\rangle _{\Sigma s}\,r^{2}$,  we get
\begin{gather}
\frac{\partial ^{2}}{\partial x^{i}\partial x^{k}}\left[ 6\left(
\frac{exp_{p}^{*}(d\Pi _{s})}{d\Pi _{Eucl.}}
 \right)+\frac{1%
}{3}\left\langle R\right\rangle _{\Sigma s}\,r^{2}\right] \underset{LGC}{=}
\\
=\frac{1}{3}\left\langle R\right\rangle _{\Sigma _{s}}\delta _{ik}-R_{ik}(p)-%
\frac{1}{2}\,\nabla _{j}\,R_{ik}(p)\,x^{j}+O(r^{2}),  \notag
\end{gather}
and by tracing
\begin{equation}
\left. \Delta _{Euc}\left[ 6\left(
\frac{exp_{p}^{*}(d\Pi _{s})}{d\Pi _{Eucl.}}
 \right)+\frac{1}{3}\left\langle
R(s)\right\rangle _{\Sigma s}\,r^{2}\right] \right| _{p}\underset{LGC}{=}%
\left\langle R\right\rangle _{\Sigma _{s}}-R(p),
\end{equation}
where $\Delta _{Euc}$ denotes the standard Euclidean Laplacian. Hence, in
local geodesic coordinates we can write
\begin{gather}
\Phi _{s}(x^{h})\underset{LGC}{=}
\left[ 6\left(
\frac{exp_{p}^{*}(d\Pi _{s})}{d\Pi _{Eucl.}}
 \right)+\frac{1}{3}\left\langle
R(s)\right\rangle _{\Sigma s}\,d(p,x)^{2}\right]= \label{lgcHess0}\\
=6-\left( R_{ik}(p)-\frac{1}{3}\left\langle
R(s)\right\rangle _{\Sigma s}\delta _{ik}\right) x^{i}x^{k}-\frac{1}{2}%
\,\nabla _{j}\,R_{ik}(p)\,x^{i}x^{k}x^{j}+O(r^{4}).
\notag
\end{gather}
Note that from this latter asymptotics we can formally compute
\begin{equation}
(Hess\Phi _{s})_{ik}(x^{h})\underset{LGC}{=}-R_{ik}(p)+\frac{1}{3}\left\langle
R(s)\right\rangle _{\Sigma s}\delta _{ik} +O(r),
\label{lgcHess}
\end{equation}
which shows that, around any given point $p$, the convexity properties of $\Phi_{s}$ are related to the sign of the
Ricci curvature. In particular, we have the
\begin{lemma}
Let $\Phi _{s}$ be a smooth solution of \ (\ref{ellipticR}) along a Ricci flow with uniformly bounded curvature operator on $\Sigma \times [0,\beta ^{*}]$,  then
\begin{equation}
Hess\,\,\Phi _{s}\geq -Ric(g(s))+\frac{1}{3}\left\langle
R(s)\right\rangle _{\Sigma s}\, g(s)
\end{equation}
in the barrier sense.
\label{normali}
\end{lemma}
\begin{proof}
Let $B_{s}(p,r)\subset (\Sigma _{s},g(s))$ be a geodesic ball of radius $r$\
centered at the generic point $p\in $\ $\Sigma _{s}$. \ For any $%
\varepsilon >0$, \ let us define
\begin{gather}
\Phi _{s}^{(\varepsilon )}(u,d(p,x))\doteq  \\
=6-\left[ Ric_{s}(u,u)|_{p}-\frac{1}{3}\left\langle R\right\rangle _{\Sigma
_{s}}+\varepsilon \,\left( \left\langle R^{2}\right\rangle _{\Sigma
_{s}}-\left\langle R\right\rangle _{\Sigma _{s}}^{2}\right) ^{\frac{1}{2}}%
\right] \,d^{2}(p,x).  \notag
\end{gather}
Away from the cut locus, the function $\Phi _{s}^{(\varepsilon )}(u,d(p,x))$ is smooth
and such that $\Phi _{s}^{(\varepsilon )}(u,d(p,x))|_{p}=\Phi _{s}(p)$. From the asymptotic expression \ of \ $\Phi _{s}(x^{h})$ in $%
B_{s}(p,r)$\ \ \ we compute
\begin{gather}
\Phi _{s}(x^{h})- \Phi
_{s}^{(\varepsilon )}(u,r)\geq  \label{utaylor}\\
-\frac{1}{2}\,|\nabla _{u}\,Ric_{s}(u,u)|_{p}|\,r^{3} +\varepsilon \,\left( \left\langle R^{2}\right\rangle _{\Sigma
_{s}}-\left\langle R\right\rangle _{\Sigma _{s}}^{2}\right) ^{\frac{1}{2}} \,r^{2}  +O(r^{4}). \notag
\end{gather}
Assume that $\nabla _{u}\,Ric_{s}(u,u)|_{p}\not=0$, (otherwise move to the next non-vanishing higher order term in the asymptotics). Along a smooth Ricci flow with uniformly bounded curvature operator on $\Sigma \times [0,\beta ^{*}]$, the derivatives $|\nabla^{(k)} \, Ric(g)|$,\, $k\geq 1$, are bounded (\cite{12}, Th. 13.1). Thus, we can assume that the $O(r^{4})$ terms in (\ref{utaylor}) are uniform in $r$, and we can define
\begin{equation}
0<\,\,r(\epsilon )\doteq 2\epsilon \,\,\frac{\left( \left\langle R^{2}\right\rangle _{\Sigma
_{s}}-\left\langle R\right\rangle _{\Sigma _{s}}^{2}\right) ^{\frac{1}{2}}}{|\nabla _{u}\,Ric_{s}(u,u)|_{p}|},
\end{equation} 
to the effect that 
\begin{equation}
\Phi _{s}(x^{h})\geq \Phi
_{s}^{(\varepsilon )}(u,r)
\end{equation}
for $0<\,r\,<r(\epsilon )$.
Finally, 
\begin{gather}
Hess\,\Phi _{s}^{(\varepsilon )}(u,d(p,x))|_{p}=-Ric_{s}(u,u)|_{p}+\frac{1}{3%
}\left\langle R(s)\right\rangle _{\Sigma _{s}}g_{s}(u,u)- \\
-\varepsilon \,\left( \left\langle R^{2}\right\rangle _{\Sigma
_{s}}-\left\langle R\right\rangle _{\Sigma _{s}}^{2}\right) ^{\frac{1}{2}%
}g_{s}(u,u).  \notag
\end{gather}
Thus, for $r$ sufficiently small, $\Phi _{s}^{(\varepsilon )}(u,d(p,x))$ is a lower barrier function for 
$\Phi _{s}(x)$. Since the base point $p\in $\ $\Sigma _{s}$ is arbitrary,  it
follows that  $Hess\,\Phi _{s}(x)\geq -Ric(s)+\frac{1}{3}\left\langle
R(s)\right\rangle _{\Sigma _{s}}g(s)$ in the barrier sense.
\end{proof}
\noindent Note that if the Ricci curvature of $(\Sigma _{\beta },g(\beta ))$ is bounded below, \emph{%
i.e.}, if there is a $K_{\beta }\in \mathbb{R}$ such that\ \ $%
R_{ik}(\beta )u^{i}u^{k}\geq K_{\beta }\,g_{ik}u^{i}u^{k}$, $\forall $ $u:\Sigma
_{\beta }\rightarrow T\Sigma _{\beta }$, then the above lemma implies that
\begin{equation}
2R_{ik}-\frac{1}{3}\left\langle
R(\beta )\right\rangle _{\Sigma \beta }g _{ik}+ (Hess\Phi _{\beta })_{ik}\geq K_{\beta }\,g_{ik}
\label{Rhessian2}
\end{equation}
a relation that will be useful in discussing the approach to equilibrium for the Fokker-Planck dynamics associated with the Ricci flow.

\subsection{Ricci flow evolution of probability measures} There is a useful consequence of the above parametrization of the curvature fluctuations which immediately shows why Fokker-Planck diffusion is natural when we
consider the evolution of a probability measure along the fiducial $d\Pi _{t}$.
\begin{lemma}
\label{L1}
For any curve of probability measures $(0,\beta _{*})\ni t\mapsto d\Omega_{t} \in Prob(\Sigma )$, absolutely
continuous with respect to $d\Pi _{t}$, the following identity holds along
the the backward volume-normalized Ricci flow  
\begin{eqnarray}
&&\frac{d }{d t}S\left[ d\Omega_{t} \parallel d\Pi _{t}\right] =-I%
\left[ d\Omega_{t} \parallel d\Pi _{t}\right] +  \label{entrdue} \\
\notag\\
&&+\int_{\Sigma _{t}}\ln \frac{d\Omega_{t} }{d\Pi _{t}}\left[ \frac{\partial }{%
\partial t}d\Omega_{t} -\Delta _{g(t)}\,d\Omega_{t} +\nabla ^{i}\left( d\Omega_{t}
\nabla _{i}\Phi _{t}\right) \right] ,  \notag
\end{eqnarray}
where 
\begin{equation}
S\left[ d\Omega_{t} \parallel d\Pi _{t}\right] \doteq \int_{\Sigma _{t}}d\Omega_{t}
\,\ln \frac{d\Omega_{t} }{d\Pi _{t}},  \label{Ltre}
\end{equation}
and
\begin{equation}
I\left[ d\Omega_{t} \parallel d\Pi _{t}\right] \doteq \int_{\Sigma _{t}}d\Omega_{t}
\,\,\nabla ^{i}\ln \frac{d\Omega_{t} }{d\Pi _{t}}\,\nabla _{i}\ln \frac{d\Omega_{t} 
}{d\Pi _{t}}
\end{equation}
respectively denote the relative entropy of $d\Omega_{t} $ with respect to $d\Pi
_{t}$ and the associated entropy generating functional. Moreover, one computes
\begin{eqnarray}
&&\frac{d }{d t}I\left[ d\Omega_{t} \parallel d\Pi _{t}\right] = \label{complic}  \\
\notag\\
&&=\int_{\Sigma _{t}}d\Omega_{t} \left[ \left| \nabla \,\ln \frac{d\Omega_{t} }{d\Pi
_{t}}\right| ^{2}+2\nabla ^{i}\ln \frac{d\Omega_{t} }{d\Pi _{t}}\nabla _{i}%
\right] \left\{ \frac{\partial }{\partial t}\ln \frac{d\Omega_{t} }{d\Pi _{t}}%
-\right.   \notag \\
\notag\\
&&\left. -\Delta _{g(t)}\ln \frac{d\Omega_{t} }{d\Pi _{t}}+\left| \nabla \ln \frac{%
d\Omega_{t} }{d\Pi _{t}}\right| ^{2}-\nabla ^{i}\Phi _{t}\nabla _{i}\ln \frac{%
d\Omega_{t} }{d\Pi _{t}}\right\} +  \notag \\
\notag\\
&&+\int_{\Sigma _{t}}d\Omega_{t} \left[ -2R^{ik}(t)+\frac{2}{3}\left\langle
R(t)\right\rangle _{\Sigma _{t}}g^{ik}(t)\right] \nabla _{i}\ln \frac{%
d\Omega_{t} }{d\Pi _{t}}\,\nabla _{k}\ln \frac{d\Omega_{t} }{d\Pi _{t}}-  \notag \\
\notag\\
&&-2\int_{\Sigma _{t}}d\Omega_{t} \left[ R^{ik}(t)+\left( Hess\,\Phi _{t}\right)
^{ik}\right] \nabla _{i}\ln \frac{d\Omega_{t} }{d\Pi _{t}}\,\nabla _{k}\ln \frac{%
d\Omega_{t} }{d\Pi _{t}}-  \notag \\
\notag\\
&&-2\int_{\Sigma _{t}}d\Omega_{t} \left| Hess\,\left( \ln \frac{d\Omega_{t} }{d\Pi _{t}%
}\right) \right| ^{2}.  \notag
\end{eqnarray}
\end{lemma}

\begin{proof}
Since by hypothesis $d\Omega_{t} $ is a probability measure absolutely continuous with respect to $d\Pi _{t}$, its total mass is
preserved along the (backward) volume-normalzed Ricci flow. By factorizing $d\Omega_{t} $ \ in terms of $d\Pi _{t}$
and by exploiting Otto's parametrization (\ref{ellipticR}), we get 
\begin{gather}
0=\frac{d }{d t}\int_{\Sigma _{t}}d\Omega_{t} =\int_{\Sigma
_{t}}d\Pi _{t}\frac{\partial }{\partial t}\left( \frac{d\Omega_{t} }{d\Pi _{t}}%
\right) +\int_{\Sigma _{t}}d\Pi _{t}\,\nabla _{i}\Phi _{t}\nabla ^{i}\left( 
\frac{d\Omega_{t} }{d\Pi _{t}}\right) =  \label{conto} \\
=\int_{\Sigma _{t}}d\Omega_{t} \,\frac{\partial }{\partial t}\ln \left( \frac{%
d\Omega_{t} }{d\Pi _{t}}\right) +\int_{\Sigma _{t}}d\Omega_{t} \,\nabla _{i}\Phi
_{t}\nabla ^{i}\ln \left( \frac{d\Omega_{t} }{d\Pi _{t}}\right) =  \notag \\
=\frac{d }{d t}\int_{\Sigma _{t}}d\Omega_{t} \,\ln \left( \frac{%
d\Omega_{t} }{d\Pi _{t}}\right) -\int_{\Sigma _{t}}\ln \left( \frac{d\Omega_{t} }{%
d\Pi _{t}}\right) \frac{\partial }{\partial t}\left( d\Omega_{t} \,\right)  
\notag \\
+\int_{\Sigma _{t}}d\Omega_{t} \,\nabla _{i}\Phi _{t}\nabla ^{i}\ln \left( \frac{%
d\Omega_{t} }{d\Pi _{t}}\right) .  \notag
\end{gather}
Integration by parts provides the identity 
\begin{gather}
\int_{\Sigma _{t}}\ln \frac{d\Omega_{t} }{d\Pi _{t}}\left[ \Delta
_{g(t)}\,d\Omega_{t} -\nabla ^{i}\left( d\Omega_{t} \nabla _{i}\Phi _{t}\right) %
\right] =  \label{Fifokker} \\
-\int_{\Sigma _{t}}d\Omega_{t} \,\,\nabla ^{i}\ln \frac{d\Omega_{t} }{d\Pi _{t}}%
\,\nabla _{i}\ln \frac{d\Omega_{t} }{d\Pi _{t}}+\int_{\Sigma _{t}}d\Omega_{t}
\,\nabla _{i}\Phi _{t}\nabla ^{i}\ln \left( \frac{d\Omega_{t} }{d\Pi _{t}}%
\right) ,  \notag
\end{gather}
which we rearrange as 
\begin{gather}
\int_{\Sigma _{t}}d\Omega_{t} \,\nabla _{i}\Phi _{t}\nabla ^{i}\ln \left( \frac{%
d\Omega_{t} }{d\Pi _{t}}\right) =I\left[ d\Omega_{t} \parallel d\Pi _{t}\right] +
\label{Fident} \\
\int_{\Sigma _{t}}\ln \frac{d\Omega_{t} }{d\Pi _{t}}\left[ \Delta
_{g(t)}\,d\Omega_{t} -\nabla ^{i}\left( d\Omega_{t} \nabla _{i}\Phi _{t}\right) %
\right] .  \notag
\end{gather}
By inserting (\ref{Fident}) in (\ref{conto}) we get(\ref{entrdue}). 
The proof of \ (\ref{complic}) is a lengthy routine computation which can be performed
along the following steps:
\begin{gather}
\frac{d }{d t}I\left[ d\Omega_{t} \parallel d\Pi _{t}\right] =%
\frac{d }{d t}\int_{\Sigma _{t}}d\Omega_{t} \,g^{ik}(t)\nabla
_{i}\ln \frac{d\Omega_{t} }{d\Pi _{t}}\,\nabla _{k}\ln \frac{d\Omega_{t} }{d\Pi _{t}}%
= \\
=\frac{d }{d t}\int_{\Sigma _{t}}d\Pi _{t}\,\left( \frac{%
d\Omega_{t} }{d\Pi _{t}}\right) ^{-1}g^{ik}(t)\nabla _{i}\frac{d\Omega_{t} }{d\Pi
_{t}}\,\nabla _{k}\frac{d\Omega_{t} }{d\Pi _{t}}=  \notag \\
=\int_{\Sigma _{t}}d\Omega_{t} \left| \nabla \,\ln \frac{d\Omega_{t} }{d\Pi _{t}}%
\right| ^{2}\left\{ R(t)-\left\langle R(t)\right\rangle _{\Sigma _{t}}+\frac{%
\partial }{\partial t}\ln \frac{d\Omega_{t} }{d\Pi _{t}}\right\} +  \notag \\
+\int_{\Sigma _{t}}d\Omega_{t} \left[ -2R^{ik}(t)+\frac{2}{3}\left\langle
R(t)\right\rangle _{\Sigma _{t}}g^{ik}(t)\right] \nabla _{i}\ln \frac{%
d\Omega_{t} }{d\Pi _{t}}\,\nabla _{k}\ln \frac{d\Omega_{t} }{d\Pi _{t}}+  \notag \\
+2\int_{\Sigma _{t}}d\Omega_{t} \,g^{ik}(t)\nabla _{i}\ln \frac{d\Omega_{t} }{d\Pi
_{t}}\,\nabla _{k}\frac{\partial }{\partial t}\ln \frac{d\Omega_{t} }{d\Pi _{t}}.
\notag
\end{gather}
The structure of \ (\ref{entrdue}) \ suggests to add and subtract, to the $\frac{\partial 
}{\partial t}\ln \frac{d\Omega_{t} }{d\Pi _{t}}$ terms in the above expression,
the quantity  
\begin{equation}
\Delta _{g(t)}\ln \frac{d\Omega_{t} }{d\Pi _{t}}+\left| \nabla \,\ln \frac{%
d\Omega_{t} }{d\Pi _{t}}\right| ^{2}-\nabla ^{i}\Phi _{t}\nabla _{i}\ln \frac{%
d\Omega_{t} }{d\Pi _{t}},
\end{equation}
which is the rewriting in terms of $\ln \frac{d\Omega_{t} }{d\Pi _{t}}$ of the
generator of the Fokker-Planck operator 
\begin{equation}
\frac{\partial }{\partial t}d\Omega_{t} -\Delta _{g(t)}d\Omega_{t} +\nabla
^{i}\left( d\Omega_{t} \nabla _{i}\Phi _{t}\right) ,
\end{equation}
appearing in (\ref{entrdue}). Applying the Bochner-Weitzenb\"{o}ck formula  
\begin{gather}
2\nabla ^{i}\ln \frac{d\Omega_{t} }{d\Pi _{t}}\nabla _{i}\Delta _{g(t)}\ln \frac{%
d\Omega_{t} }{d\Pi _{t}}=\Delta _{g(t)}\left| \nabla \ln \frac{d\Omega_{t} }{d\Pi
_{t}}\right| ^{2}- \\
-2\left| Hess\,\ln \frac{d\Omega_{t} }{d\Pi _{t}}\right| ^{2}-2R^{ik}(t)\,\nabla
_{i}\ln \frac{d\Omega_{t} }{d\Pi _{t}}\nabla _{k}\ln \frac{d\Omega_{t} }{d\Pi _{t}},
\notag
\end{gather}
and the identities 
\begin{gather}
\int_{\Sigma _{t}}d\Omega_{t} \left| \nabla \,\ln \frac{d\Omega_{t} }{d\Pi _{t}}%
\right| ^{2}\Delta _{g(t)}\ln \frac{d\Omega_{t} }{d\Pi _{t}}+\int_{\Sigma
_{t}}d\Omega_{t} \left| \nabla \,\ln \frac{d\Omega_{t} }{d\Pi _{t}}\right| ^{4}+ \\
+\int_{\Sigma _{t}}d\Omega_{t} \,\Delta _{g(t)}\left| \nabla \,\ln \frac{d\Omega_{t} 
}{d\Pi _{t}}\right| ^{2}+2\int_{\Sigma _{t}}d\Omega_{t} \,\nabla ^{i}\ln \frac{%
d\Omega_{t} }{d\Pi _{t}}\nabla _{i}\left| \nabla \,\ln \frac{d\Omega_{t} }{d\Pi _{t}}%
\right| ^{2}=0,  \notag
\end{gather}
\begin{gather}
\nabla ^{i}\ln \frac{d\Omega_{t} }{d\Pi _{t}}\nabla _{i}\left( \nabla ^{k}\Phi
_{t}\nabla _{k}\ln \frac{d\Omega_{t} }{d\Pi _{t}}\right) -\frac{1}{2}\nabla
^{k}\Phi _{t}\nabla _{k}\left| \nabla \,\ln \frac{d\Omega_{t} }{d\Pi _{t}}%
\right| ^{2}= \\
=\nabla ^{i}\ln \frac{d\Omega_{t} }{d\Pi _{t}}\nabla ^{k}\ln \frac{d\Omega_{t} }{%
d\Pi _{t}}\left( Hess\,\Phi _{t}\right) _{ik},  \notag
\end{gather}
(the former obtained by iterated integrations by parts and the second by
direct computation \cite{17}), one eventually gets the stated result.
\end{proof}

\subsection{Fokker-Planck diffusion along the Ricci flow} From Lemma \ref{L1}, we immediately get the following

\begin{theorem}
The Fokker-Planck diffusion $(d\Omega _{t})_{t\geq
0}$ generated, along the backward, volume-preserving, Ricci flow, by 
\begin{equation}
\frac{\partial }{\partial t}\left( d\Omega _{t}\right) =\Delta
_{g(t)}d\Omega _{t}-\nabla ^{i}\left( d\Omega _{t}\nabla _{i}\Phi
_{t}\right) ,
\label{FoPl}
\end{equation}
has the following properties:

\noindent (i)\; It is a gradient--like flow for the relative entropy functional $S\left[ d\Omega
_{t}\parallel d\Pi _{t}\right] $, i.e.,
\begin{equation}
\frac{d }{d t}\int_{\Sigma _{t}}d\Omega _{t}\,\ln \frac{%
d\Omega _{t}}{d\Pi _{t}}=-\int_{\Sigma _{t}}d\,\Omega _{t}\nabla _{i}\ln
\left( \frac{d\Omega _{t}}{d\Pi _{t}}\right) \nabla ^{i}\ln \left( \frac{%
d\Omega _{t}}{d\Pi _{t}}\right)\;.
\label{gradS}
\end{equation}

\noindent (ii)\;The corresponding evolution for the Radon--Nikodym derivative $\frac{d\Omega _{t}}{d\Pi _{t}}$
\begin{equation}
\frac{\partial }{\partial t}\frac{d\Omega _{t}}{d\Pi _{t}}=\Delta _{g(t)}%
\frac{d\Omega _{t}}{d\Pi _{t}}-\nabla ^{i}\Phi _{t}\nabla _{i}\frac{d\Omega
_{t}}{d\Pi _{t}},
\label{FPdensity}
\end{equation}
is the gradient flow, (in the $L^{2}(\Sigma ,d\Pi _{t})$ sense), of $S\left[ d\Omega
_{t}\parallel d\Pi _{t}\right] $ with respect to the inner product $\langle\ldots,\ldots\rangle_{d\Omega _{t}\backslash d\Pi _{t}}$, defined by (\ref{inner}). 

\noindent (iii)\;The associated entropy production functional $I\left[ d\Omega _{t}\parallel d\Pi _{t}\right]$  satisfies the 
differential inequality
\begin{equation}
\frac{d }{d t}I\left[ d\Omega _{t}\parallel d\Pi _{t}\right] 
\leq -2\,K_{t}\,I\left[ d\Omega _{t}\parallel d\Pi _{t}\right] 
\label{ineqI}
\end{equation}
where $K_{t}\,\in\mathbb{R}$ is the ($t$-dependent) lower bound of the Ricci curvature.

\noindent (iv)\; If along the given  Ricci flow $\beta \mapsto g_{ab}(\beta )$ the diameter stays uniformly bounded, i.e.,  $\sup_{\beta \geq 0}\;\{diam\,(\Sigma ,g(\beta ))\}\,\doteq \,\overline{diam}\,\,<\infty $, then  
\begin{equation}
S[d\Omega _{t}\parallel d\Pi _{t}]\geq \frac{1}{2}\left[ \frac{%
D_{2}^{W}(d\Omega _{t},d\Pi _{t})}{\overline{diam}}\right] ^{4}\;,
\label{toimprove}
\end{equation}
and the quadratic Wasserstein distance $D_{2}^{W}(d\Omega _{t},d\Pi _{t})$ is weakly monotonically decreasing along the backward Ricci flow.
\label{teorino}
\end{theorem}
\begin{proof}

\emph{(i)}\;The first part of the theorem, and in particular the gradient--like nature of $(d\Omega _{t})_{t\geq
0}$, is a direct computational consequence of Lemma \ref{L1}.

\noindent \emph{(ii)}\; It is easily verified that the absolutely continuous curve $t\mapsto d\Omega _{t}$ solution of (\ref{FoPl}) induces the evolution (\ref{FPdensity}) on the associated Radon-Nikodym derivatives $(\frac{d\Omega _{t}}{d\Pi  _{t}})_{t\geq 0}$. The proof that this  is  a gradient flow with respect to the inner product $\langle\ldots,\ldots\rangle_{d\Omega _{t}\backslash d\Pi _{t}}$ follows from the analysis in section \ref{Ottop}, (see eq. (\ref{GFS})). In particular, if we denote by $\Psi _{t}\in C_{b}(\Sigma ,\mathbb{R})\slash\mathbb{R}$ the tangent vector to $(d\Omega _{t})_{t\geq
0}$, \emph{i.e.}, the solution, for each given $t=s$, of the elliptic PDE
\begin{equation}
\left.\frac{\partial }{\partial t}\,d\Omega _{t}\right|_{t=s}=\,-g^{ik}(s)\nabla _{i}\left(d\Omega _{t}\,\nabla _{k}\Psi _{s} \right)\;,
\label{tomega}
\end{equation}
then, according to (\ref{convDer}) and (\ref{RDtangent3}), the tangent vector to the curve of Radon--Nikodym derivatives $(\frac{d\Omega _{t}}{d\Pi _{t}})_{t\geq 0}$, at the given value of $t$,  is provided by $(\Psi _{t}-\Phi _{t})$. The gradient flow condition with respect to the inner product (\ref{inner}) takes the form (\ref{gradcond}), \emph{i.e.}
\begin{equation}
\left\langle (\Psi _{t }-\Phi  _{t }),\,\xi  \right\rangle_{\frac{d\Omega  _{t }}{d\Pi _{t }}}+
\left\langle Grad\, S[d\Omega _{t}\parallel d\Pi_{t}],\,\xi  \right\rangle_{\frac{d\Omega  _{t }}{d\Pi _{t }}}=0\;,
\end{equation}

\noindent $\forall \xi \in T_{\frac{d\Omega  _{t }}{d\Pi _{t }}}Prob(\Sigma_{t})$, which is equivalent, in the ${L^{2}(\Sigma ,{d\Pi _{t }})}$  sense, to (\ref{FPdensity}). Note that
by comparing (\ref{tomega}) with (\ref{FoPl}), one gets the relation $\nabla _{i}(\Psi _{t}-\Phi _{t})$ $=$ $-\nabla _{i}\ln \left( \frac{%
d\Omega _{t}}{d\Pi _{t}}\right)$. Thus, we can equivalently rewrite the entropy production functional as the Otto norm of the vector $(\Psi _{t}-\Phi _{t})$, \emph{i.e.}, 
\begin{equation}
I\left[ d\Omega _{t}\parallel d\Pi _{t}\right]=\int_{\Sigma _{t}}\left|\nabla\,(\Psi _{t}-\Phi _{t}) \right|^{2}\,
\frac{d\Omega _{t}}{d\Pi  _{t}}\,d\Pi  _{t}\;.
\label{lntangent}
\end{equation}

 \noindent (iii)\;From Lemma \ref{L1} we have  
\begin{gather}
\frac{d }{d t}I\left[ d\Omega _{t}\parallel d\Pi _{t}\right] = \label{gradI}
\\
=\int_{\Sigma _{t}}d\Omega _{t}\left[ -2R^{ik}(t)+\frac{2}{3}\left\langle
R(t)\right\rangle _{\Sigma _{t}}g^{ik}(t)\right] \nabla _{i}\ln \frac{%
d\Omega _{t}}{d\Pi _{t}}\,\nabla _{k}\ln \frac{d\Omega _{t}}{d\Pi _{t}}- 
\notag \\
-2\int_{\Sigma _{t}}d\Omega _{t}\left[ R^{ik}(t)+\left( Hess\,\Phi
_{t}\right) ^{ik}\right] \nabla _{i}\ln \frac{d\Omega _{t}}{d\Pi _{t}}%
\,\nabla _{k}\ln \frac{d\Omega _{t}}{d\Pi _{t}}-  \notag \\
-2\int_{\Sigma _{t}}d\Omega _{t}\left| Hess\,\left( \ln \frac{d\Omega _{t}}{%
d\Pi _{t}}\right) \right| ^{2},  \notag
\end{gather}
which, according to Lemma \ref{normali} and (\ref{Rhessian2}), directly yields (\ref{ineqI}).

\noindent (iv)\; As long as the diameter $diam\,(\Sigma ,g(\beta ))$ of $(\Sigma ,g(\beta ))$ remains uniformly bounded along the given Ricci flow, we have the Talagrand--like inequality (\ref{Pinsk2}) from which (\ref{toimprove}) immediately follows.
Note that for a Ricci flow with uniformly bounded Ricci curvature $|Ric(\beta )|\leq M$ we have the elementary  bound for $diam\,(\Sigma ,g(\beta ))$, (see \cite{21aa} and also \cite{12}),
\begin{equation}
diam\,(\Sigma ,g(0))\,e^{-2M\,\beta }\leq\,diam\,(\Sigma ,g(\beta ))\,\leq diam\,(\Sigma ,g(0))\,e^{2M\,\beta }\;, 
\end{equation}
(the factor 2 is due to the volume normalization of the flow). Recently, P. Topping \cite{21aa} has obtained an improved control on $diam\,(\Sigma ,g(\beta ))$ in terms of suitable averages of the scalar curvature. 
\end{proof}

\noindent In its simplest form,  the rate of convergence of a solution of 
\begin{equation}
\left\{ 
\begin{tabular}{l}
$\frac{\partial }{\partial t}\frac{d\Omega _{t}}{d\Pi _{t}}=\Delta _{g(t)}%
\frac{d\Omega _{t}}{d\Pi _{t}}-\nabla ^{i}\Phi _{t}\nabla _{i}\frac{d\Omega
_{t}}{d\Pi _{t}},$ \\ 
\\ 
$\frac{d\Omega _{t}}{d\Pi _{t}}|_{t=0}=\frac{d\Omega _{0}}{d\Pi _{0}},$%
\end{tabular}
\right. 
\end{equation}
to the stationary state $d\Pi _{t}$ is governed  by a curvature condition which is naturally suggested by the structure of equation (\ref{gradI}), and which, according to lemma \ref{normali}, \, is equivalent to the  positivity of the Ricci tensor. Note that, according to the characterization (\ref{ellipticR}) of the potential $\Phi _{t}$,  also $(d\Pi  _{t})_{t\geq 0}$  solves the Fokker--Planck equation (\ref{FoPl}), (with the initial datum $d\Pi  _{0}$), along the backward Ricci flow. If $(d\Omega _{t})_{t\geq 0}$ is a solution of (\ref{FoPl}) then we call $(d\Omega _{t},\,\,d\Pi  _{t} )_{t\geq 0}$ a conjugated Fokker-Planck pair along the backward Ricci flow.   Taking into account this elementary observation,  we get the following result 
\begin{lemma}
Let $\beta \rightarrow g_{ab}(\beta )$,  $%
\beta \in \lbrack 0,\infty )$, a given Ricci flow metric starting on a manifold $(\Sigma ,g)$ of positive Ricci curvature, and let $(d\Omega _{t},\,\,d\Pi  _{t} )_{t\geq 0}$ be the conjugated Fokker--Planck pair, solution of the Fokker-Planck equation (\ref{FoPl}). 
Then, along the backward Ricci flow,  
the entropy functional $S\left[ d\Omega _{t}\parallel d\Pi _{t}%
\right] $ decreases exponentially fast according to
\begin{equation}
S\left[ d\Omega _{t}\parallel d\Pi _{t}\right] \leq S\left[ d\Omega
_{0}\parallel d\Pi _{0}\right] e^{-\frac{2}{3}\lambda _{\inf }\,\,t},
\label{Sequilibrium}
\end{equation}
\noindent where  $\frac{1}{3}\lambda _{\inf }\doteq \inf_{\beta\geq 0}\,\,\left\{K_{\beta }>0:\,Ric(\beta )\geq K_{\beta }\,g(\beta ) \right\}$. Moreover, the Talagrand inequality
\begin{equation}
S\left[ d\Omega _{t}\parallel d\Pi _{t}\right]\,\geq \frac{\lambda _{\inf }}{6}\,\left[D_{2}^{W}(d\Omega _{t},d\Pi _{t})\right]^{2}\;,
\label{talagrand} 
\end{equation}
holds, and $S\left[ d\Omega _{t}\parallel d\Pi _{t}\right]$ is a convex function along the backward Ricci flow. 
\end{lemma}
\begin{proof}
Since the positivity condition on the Ricci tensor is preserved by the Ricci flow and yields long time existence \cite{11}, 
we have  $Ric(\beta )\geq K_{\beta }\,g(\beta )$, 
$K_{\beta }\geq c>0$ along the flow. According to (\ref{Rhessian2}), 
such positivity implies the condition
\begin{equation}
2Ric(\beta )-\frac{1}{3}\left\langle R(\beta )\right\rangle _{\Sigma _{\beta
}}g(\beta )+Hess\,\Phi _{\beta }\geq \frac{\lambda _{inf }}{3}\,g(\beta ),
\label{noBE}
\end{equation}
on the curvature terms entering (\ref{gradI}). It follows that
\begin{equation}
\frac{d }{d t}I\left[ d\Omega _{t}\parallel d\Pi _{t}\right]
\leq -\frac{2}{3}\lambda _{\inf }\,I\left[ d\Omega _{t}\parallel d\Pi _{t}\right] ,
\label{partialIa}
\end{equation}
which implies that the entropy dissipation functional  $I\left[ d\Omega
_{t}\parallel d\Pi _{t}\right] $ decreases exponentially fast according to 
\begin{equation}
I\left[ d\Omega _{s}\parallel d\Pi _{s}\right] \leq I\left[ d\Omega
_{t}\parallel d\Pi _{t}\right] \,e^{-\frac{2}{3}\lambda _{\inf }\,\,(s-t)},
\label{Iomega}
\end{equation}
for any $(s-t)\geq 0$. Let $t_{eq}$ (this may be finite or
infinite) the value of $t$ for which $d\Omega _{t}$ attains equilibrium
(i.e. $d\Omega _{t}=d\Pi _{t}$ for $t\geq t_{eq}$).\ If we integrate (\ref
{Iomega}) over $s$\ from $t$ to $t_{eq}$, 
\begin{equation}
\int_{t}^{t_{eq}}I\left[ d\Omega _{s}\parallel d\Pi _{s}\right] \,ds\leq 
\frac{3}{2\,\lambda _{\inf }}\left( 1-e^{-\frac{2}{3}\lambda _{\inf
}\,\,(t_{eq}-t)}\right) \,I\left[ d\Omega _{t}\parallel d\Pi _{t}\right] .
\end{equation}
and take into account that  $I\left[ d\Omega _{s}\parallel d\Pi _{s}\right]
=-\frac{d }{d s}S\left[ d\Omega _{s}\parallel d\Pi _{s}\right] 
$, together with $S\left[ d\Omega _{t_{eq}}\parallel d\Pi _{te_{q}}\right] \doteq 0$,\ and 
$\left( 1-e^{-\frac{2}{3}\lambda _{\inf }\,\,(t_{eq}-t)}\right) \leq 1$, we get
\begin{equation}
S\left[ d\Omega _{t}\parallel d\Pi _{t}\right] \leq \frac{3}{2\,\lambda
_{\inf }}\,I\left[ d\Omega _{t}\parallel d\Pi _{t}\right] .
\label{lsilambda}
\end{equation}
Since the time $t$ is arbitrary, this establishes  that a logarithmic
Sobolev inequality, of constant $\lambda _{\inf }$, holds\ for the diffusion
process $\{d\Omega _{t}\}_{t\geq 0}$.
By inserting (\ref{lsilambda})\ \ in (\ref{gradS})
we immediatey get the exponential trend to equilibrium (\ref{Sequilibrium}). Since $(d\Omega _{t},\,\,d\Pi  _{t} )_{t\geq 0}$ are a conjugated Fokker--Planck pair, according to \cite{17}, the validity of the logarithmic Sobolev inequality on $(\Sigma, g(t) )$ implies the Talagrand inequality
and hence (\ref{talagrand}) follows. The convexity of $S\left[ d\Omega _{t}\parallel d\Pi _{t}\right]$ is a direct consequence of 
(\ref{partialIa}) and (\ref{lsilambda}) above. Explicitly, from (\ref{partialIa}) and $I\left[ d\Omega _{s}\parallel d\Pi _{s}\right]
=-\frac{d }{d s}S\left[ d\Omega _{s}\parallel d\Pi _{s}\right]$, we get
\begin{equation}
\frac{d ^{2}}{d t^{2}}
S\left[ d\Omega _{t}\parallel d\Pi _{t}\right] \geq \left(\frac{2}{3}\,\lambda _{inf} \right)^{2}
S\left[ d\Omega _{t}\parallel d\Pi _{t}\right], 
\end{equation}
which by Talagrand inequality (\ref{talagrand}) yields
\begin{equation}
\frac{d ^{2}}{d t^{2}}
S\left[ d\Omega _{t}\parallel d\Pi _{t}\right] \geq \frac{1}{4}\left(\frac{2}{3}\,\lambda _{inf} \right)^{3}
\left[D_{2}^{W}(d\Omega _{t},d\Pi _{t})\right]^{2},
\label{sturmete}
\end{equation}
and $S\left[ d\Omega _{t}\parallel d\Pi _{t}\right]$ is $t$-displacement convex along the Ricci flow.
\end{proof}
Note that if we introduce the adimensional variable $\eta \doteq \frac{1}{3}\lambda _{inf}\,t$ then (\ref{sturmete}) can be equivalently rewritten as 
\begin{equation}
\frac{d ^{2}}{d \eta ^{2}}
S\left[ d\Omega _{t}\parallel d\Pi _{t}\right] \geq \frac{2}{3}\lambda _{inf} 
\left[D_{2}^{W}(d\Omega _{t},d\Pi _{t})\right]^{2},
\end{equation} 
which is equivalent to Sturm's $K$-convexity \cite{23} of $S\left[ d\Omega _{t}\parallel d\Pi _{t}\right]$, (for $K=\frac{1}{3}\,\lambda _{inf}$), (see (\ref{KK})). The point here is that, typically, the $K$-convexity of a relative entropy functional  holds along the Wasserstein geodesics of the metrized probability space $(Prob(\Sigma,g),\,D_{2}^{W}(\;,\;))$. In particular, it was conjectured to hold for Riemannian manifolds with non-negative Ricci curvature by F. Otto and C. Villani \cite{17}, (in the case $K=0$, the conjecture has been proven in \cite{6a} whereas in the general case ( for any $K\in\mathbb{R}$) in \cite{23}). These remarks suggest that the Fokker-Planck diffusion $\{d\Omega _{t}\}_{t\geq 0}$ is strictly connected with Wasserstein geodesics in the bundle $Prob(\Sigma )$. To discuss to what extent this is the case, let us recall that, at any given $t$, the tangent vector to the curve of Radon--Nikodym derivatives $t\mapsto \frac{d\Omega _{t}}{d\Pi _{t}}$, is provided by $-\ln\left(\frac{d\Omega _{t}}{d\Pi _{t}} \right)$, (see (\ref{lntangent})). In analogy with the characterization of the parameter $\tau$ characterizing Perelman's diffusion $(d\varpi (t))_{t\geq 0}$, let $\mu _{0}$ denote the (squared) length scale over which the probability measure $d\Omega _{t=0}$ is concentrated. Roughly speaking $(4\pi \mu _{0})^{\frac{3}{2}}$ is the typical volume of a  domain $B\subset$ whose $d\Omega _{0}$-measure is not exponentially small and yields for full measure if slightly blown up. Let us consider the adimensional parameter
\begin{equation}
\varepsilon \doteq \frac{\lambda _{inf} }{3}\,(4\pi \mu _{0})^{\frac{2}{3}}\;,
\end{equation}
and let us  $\epsilon $--rescale the vector $-\ln\left(\frac{d\Omega _{t}}{d\Pi _{t}} \right)$ according to
\begin{equation}
u_{\varepsilon }(\varepsilon \,t)\doteq -2\,\varepsilon \,\ln\left(\frac{d\Omega _{t}}{d\Pi _{t}} \right)\;.
\label{epsilocal}
\end{equation}

\noindent From equation (\ref{FPdensity}) it easily follows that $(u_{\varepsilon }(\varepsilon \,t))_{t\geq 0}$ evolves according to the viscous Hamilton--Jacobi equation \cite{22, 22a}
\begin{equation}
\frac{\partial \,u_{\varepsilon } }{\partial t}\,+\frac{\left|u_{\varepsilon } \right|^{2}}{2}=\,\varepsilon \,\left(\Delta _{g(t)}%
u_{\varepsilon }-\nabla ^{i}\Phi _{t}\nabla _{i}u_{\varepsilon }\right)\;,
\label{HJFP}
\end{equation}
with (a smooth) initial datum $u_{\varepsilon }(x,t=0)=U_{\varepsilon }(x)$, $x\in \Sigma $.
When the parameter $\varepsilon $ defined by (\ref{epsilocal}) is small, one may discuss the solution of (\ref{HJFP}) by the so--called vanishing viscosity method. Qualitatively, this implies that when $\varepsilon \rightarrow 0$, the rescaled vector $u_{\varepsilon }$ approaches the Hopf--Lax solution 
\begin{equation}
\inf_{y\in \Sigma }\left[U_{\varepsilon }(y)+\frac{1}{2\,t}\,d_{t}(x,y)^{2} \right]\;,\;t>0,\;\;x\in \Sigma\;, 
\label{hopf}
\end{equation}
of the Hamilton-Jacobi equation
\begin{equation}
\left\{ 
\begin{tabular}{l}
$\frac{\partial \,u_{\varepsilon } }{\partial t}\,+\frac{\left|u_{\varepsilon } \right|^{2}}{2}=0$ \\ 
\\ 
$u_{\varepsilon }(x,t=0)=U_{\varepsilon }(x).$%
\end{tabular}
\right.   \label{trueHJ}
\end{equation}
This is admittedly rather vague since, in our setting, the distance $d_{t}(x,y)$ varies with $t$, along the backward Ricci flow, and the viscosity solutions must take this dependence into account. However, for  $\varepsilon\,<<1 $, \emph{i.e.} if the probability measure $d\Omega _{t}$ is initially concentrated on a set which is small with respect to the radius of curvature of $(\Sigma ,g(t=0))$, the Fokker--Planck diffusion $(d\Omega _{t})_{t\geq 0}$, behaves, for $t$ sufficiently small, as if occurring in the fixed probability space $Prob\,(\Sigma ,g(t=0))$, with a shadow of the Ricci flow still present through the forcing potential $\Phi _{t=0}$. In such a case,  the viscosity interpretation of (\ref{hopf}) is more justified and, according to theorem \ref{LottLength}, (see (\ref{HamJac})), one can reach the conclusion that $(d\Omega _{t})_{t\geq 0}$, for $t$ small enough, approximates a geodesic in the Wasserstein space  $(Prob(\Sigma,g(0)),\,D_{2}^{W}(\;,\;))$. From a more geometrical point of view, one is here approximating the curve $t\mapsto \frac{d\Omega _{t}}{d\Pi _{t}}$ with the push--forward of $\left.\frac{d\Omega _{t}}{d\Pi _{t}}\right|_{t=0}$ under the action of the semigroup defined by the infinitesimal generator $-\frac{1}{2}|\nabla U_{\varepsilon }|^{2}$, (where the norm is taken with respect to $(\Sigma ,g(t=0))$). It is a known fact that, whenever one moves a measure density through the push--forward action of the exponential  of a smooth function, one gets a geodesic in the appropriate Wasserstein space, (again, I wish to thank A. Figalli for useful remarks in this connection).

\noindent What is missing in such an approximation argument  is the explicit role of the backward Ricci flow. In this connection, the basic step to take  is the characterization of Wasserstein geodesics not just in a fixed probability space $(Prob(\Sigma,g),\,D_{2}^{W}(\;,\;))$, but rather in the bundle $Prob(\Sigma )$ over the space of Riemannian metrics $\mathcal{R}iem(\Sigma )$. Roughly speaking, we  expect that  the Hamilton-Jacobi condition (\ref{HamJac}) of theorem  \ref{LottLength} is an approximation to a more general geodesic equation in $Prob(\Sigma )$, more or less like straight lines approximate geodesics in Riemannian geometry.  Correspondingly, the Hopf-Lax representation will be an approximation to the exponential mapping in $Prob(\Sigma )$. These remarks are strongly supported by the fact that, as we have seen above, the natural diffusion process along the backward Ricci flow is the viscous Hamilton--Jacobi equation (\ref{HJFP}), where the viscosity parameter $\varepsilon $ is naturaly characterized by the lower bound of the Ricci curvature. Here, we see a recurring theme in the Fokker--Planck dynamics of the conjugated pair $(d\Omega _{t},\,d\Pi _{t})_{t\geq 0}$ along the (backward) Ricci flow: \emph{(i)}\;The Ricci curvature controls the geodesic convexity for the corresponding relative entropy; \emph{(ii)}\; It parametrizes the viscosity solutions of the Hamilton-Jacobi equation associated with the Fokker--Planck  diffusion; \emph{(iii)}\; It naturally affects the logarithmic Sobolev inequalities controlling the Wasserstein distance between $(d\Omega _{t},\,d\Pi _{t})_{t\geq 0}$. Points of contact between diffusion, geodesic convexity, Hamilton--Jacobi theory, and LSI are well--known (see \emph{e.g.}, \cite{22,22a} for a discussion and relevant references), however here their relation seem to come to full circle. In our opinion,  this is a serious indication of the existence of deeper connections between the geometry of optimal transport and Ricci flow theory. In particular such connections point to the possibility of adopting the geometry optimal trasportation for extending the Ricci flow to metric spaces more general that Riemannian manifolds.

\vfill\eject
\section{Appendix: Comparison with Perelman's heat flow}
It is worthwhile to compare the Fokker-Planck diffusion (\ref{FoPl}) with Perelman's flow (\ref{ourp}).
If we apply lemma \ref{L1} to $(d\varpi (t))_{t\geq 0}$ we
immediately get the following
\begin{lemma}
For a given Ricci flow metric $\beta \longmapsto g(\beta )$, $\beta \in
\lbrack 0,T)$, and for any chosen $\beta ^{\ast }\in \lbrack 0,T)$,\ let $%
t\longmapsto d\varpi (t)$, $t\doteq \beta ^{\ast }-\beta $, be a solution of
\ the parabolic PDE (\ref{heat}). Then, the relative entropy functional 
\begin{equation}
S\left[ d\varpi (t)\parallel d\Pi _{t}\right] \doteq \int_{\Sigma
_{t}}d\varpi (t)\ln \frac{d\varpi (t)}{d\Pi _{t}}
\end{equation}
\ varies along the fiducial flow $t\longmapsto d\Pi _{t}$ according to \ \ \ 
\begin{gather}
\frac{d }{d t}\int_{\Sigma _{t}}d\varpi (t)\ln \frac{d\varpi
(t)}{d\Pi _{t}}=-I\left[ d\varpi (t)\parallel d\Pi _{t}\right] -
\label{monot2} \\
-\int_{\Sigma _{t}}\nabla ^{i}\frac{d\varpi (t)}{d\Pi _{t}}\nabla _{i}\Phi
_{t}\,d\Pi _{t}.\;\;   \notag
\end{gather}
\label{gradfail}
\end{lemma}
Note in particular that $S\left[ d\varpi (t)\parallel d\Pi _{t}\right]$, as compared to $S\left[ d\Omega _{t}\parallel d\Pi _{t}\right]$,  is not weakly monotonic. The term responsible for such a lack of monotonicity,
(and of the fact that (\ref{heat}) is not the gradient flow of $S\left[ d\varpi (t)\parallel d\Pi _{t}\right]$), is the scalar product 
$\left\langle  \Phi _{t}, \frac{d\varpi (t)}{d\Pi _{t}}\right\rangle _{d\Pi _{t}}$. If we replace $d\varpi (t)$ with $d\Omega _{t}$, this is basically the drift term driving Fokker-Planck diffusion in (\ref{FoPl}). Such a term describes also  how the measure $d\varpi (t)$ localizes the fluctuations in
the scalar curvature along the Ricci flow. 
\begin{lemma}
For a given Ricci flow metric $\beta \longmapsto g(\beta )$, $\beta \in
\lbrack 0,T)$, and for any chosen $\beta ^{\ast }\in \lbrack 0,T)$,\ let $%
t\longmapsto d\varpi (t)$, $t\doteq \beta ^{\ast }-\beta $, be a solution of
\ the parabolic equation (\ref{heat}) corresponding to the initial datum $d\varpi
_{0}=(d\varpi _{0}/d\Pi _{\beta ^{\ast }})d\Pi _{\beta ^{\ast }}$. Then 
\begin{equation}
\left\langle  \Phi _{t}, \frac{d\varpi (t)}{d\Pi _{t}}%
\right\rangle _{d\Pi _{t}}=\int_{\Sigma }d\varpi (t)\left( R(t)-\left\langle
R(t)\right\rangle \right) ,  \label{omegaR}
\end{equation}
and 
\begin{gather}
\frac{d }{d t}\left\langle  \Phi _{t}, \frac{%
d\varpi (t)}{d\Pi _{t}}\right\rangle _{d\Pi _{t}}=-2\int_{\Sigma }d\varpi
(t)\left( \left| \widehat{Ric}(t)\right| ^{2}-\langle \left| \widehat{Ric}%
(t)\right| ^{2}\rangle _{\Sigma _{t}}\right)   \label{secscal} \\
-\frac{2}{3}\int_{\Sigma }d\varpi (t)R(t)\left[ R(t)-\langle R(t)\rangle
_{\Sigma _{t}}\right] -\frac{1}{3}[\langle R(t)^{2}\rangle _{\Sigma _{t}}-%
\langle R(t)\rangle _{\Sigma _{t}}^{2}].  \notag
\end{gather}
\end{lemma}

\begin{proof}
Relation (\ref{omegaR}) follows from (\ref{ellipticR}) by a straightforward
integration by parts 
\begin{gather}
\left\langle  \Phi _{t}, \frac{d\varpi (t)}{d\Pi _{t}}%
\right\rangle _{d\Pi _{t}}=\int_{\Sigma }d\Pi _{t}\nabla ^{i}\left( \frac{%
d\varpi (t)}{d\Pi _{t}}\right) \nabla _{i}\Phi _{t}=  \label{scalar1} \\
=\int_{\Sigma }\nabla ^{i}\left( d\varpi (t)\right) \nabla _{i}\Phi
_{t}=\int_{\Sigma }\nabla ^{i}\left( d\varpi (t)\nabla _{i}\Phi _{t}\right)
-\int_{\Sigma }d\varpi (t)\Delta \Phi _{t}=  \notag \\
=\int_{\Sigma }d\varpi (t)\left( R(t)-\left\langle R(t)\right\rangle \right)
.  \notag
\end{gather}
The evolution equation for scalar curvature (\ref{evscalcurv}) provides 
\begin{equation}
\frac{d }{d \beta }\langle R(\beta )\rangle _{\Sigma _{\beta
}}=2\langle \left| \widehat{Ric}(\beta )\right| ^{2}\rangle _{\Sigma _{\beta
}}-\frac{1}{3}\left( \langle R(\beta )^{2}\rangle _{\Sigma _{\beta
}}-\langle R(\beta )\rangle _{\Sigma _{\beta }}^{2}\right) .
\end{equation}
From these latter relation, (\ref{evscalcurv}) and (\ref{omegaR}) one directly computes 
\begin{gather}
\frac{d }{d t}\left\langle  \Phi _{t}, \frac{%
d\varpi (t)}{d\Pi _{t}}\right\rangle _{d\Pi _{t}}= \\
=\int_{\Sigma }R(t)\Delta (d\varpi (t))+\int_{\Sigma }d\varpi (t)\frac{%
\partial }{\partial t}R(t)-\frac{d }{d t}\left\langle
R(t)\right\rangle   \notag \\
=-2\int_{\Sigma }d\varpi (t)\left( \left| \widehat{Ric}(t)\right|
^{2}-\langle \left| \widehat{Ric}(t)\right| ^{2}\rangle _{\Sigma
_{t}}\right) -  \notag \\
-\frac{2}{3}\int_{\Sigma }d\varpi (t)R(t)\left[ R(t)-\langle R(t)\rangle
_{\Sigma _{t}}\right] +  \notag \\
-\frac{1}{3}[\langle R(t)^{2}\rangle _{\Sigma _{t}}-\langle
R(t)\rangle _{\Sigma _{t}}^{2}],  \notag
\end{gather}
which provides  the stated result (\ref{secscal}). 
\end{proof}
Note in particular that if we choose for (\ref{heat}) the
initial datum  $d\varpi (t=0)=d\Pi _{\beta ^{\ast }}$ we get, 
\begin{equation}
\left. \frac{d }{d t}\left\langle  \Phi _{t}, 
\frac{d\varpi (t)}{d\Pi _{t}}\right\rangle _{d\Pi _{t}}\right|
_{t=0}=-\left( \langle R(\beta ^{\ast })^{2}\rangle _{\Sigma _{\beta ^{\ast
}}}-\langle R(\beta ^{\ast })\rangle _{\Sigma _{\beta ^{\ast }}}^{2}\right),
\end{equation}
which again shows the role that mean square fluctuations in scalar curvature have in controlling the 
concentration mechanism of the measure $d\varpi (t)$. We can trace here the difference between the standard Perelman flow characterizing the backward diffusion of $d\varpi (t)$ and the Fokker-Planck diffusion of $d\Omega _{t}$. The former feels curvature fluctuations in a more indirect way as a forcing effect deforming the
trajectory of $d\varpi (t)$ in $Prob(\Sigma )$. Such a forcing behavior is made manifest by the fact that the evolution of $d\varpi (t)$ is not the gradient flow of the associated relative entropy, and, according to lemma \ref{gradfail}, the failure of being gradient is exactly provided by the term $\left\langle  \Phi _{t}, 
\frac{d\varpi (t)}{d\Pi _{t}}\right\rangle _{d\Pi _{t}}$. Conversely, the diffusion $(d\Omega _{t})_{t\geq 0}$ has the curvature fluctuations taken care of by turning the forcing term $\left\langle  \Phi _{t}, 
\frac{d\varpi (t)}{d\Pi _{t}}\right\rangle _{d\Pi _{t}}$ into the drift term $\left\langle  \Phi _{t}, 
\frac{d\Omega _{t}}{d\Pi _{t}}\right\rangle _{d\Pi _{t}}$   which renormalizes Perelman's $\{d\varpi (t)\}_{t\geq 0}$ into the Fokker-Planck diffusion 
$\{d\Omega _{t}\}_{t\geq 0}$.

\bigskip

\subsection*{Aknowledgements}
The author would like to thank T. Buchert and D. Glickenstein for useful conversations in the preliminary stage of preparation of this paper. Discussions with G. Bellettini, A. Figalli,  G. Savar\' e, G. Toscani, and the remarks of the referee have been extremely helpful in improving the presentation.

\bibliographystyle{amsplain}

\end{document}